\documentclass[11pt]{article}

\author{M.~Kazarian\footnote{Partially supported by RFBR grant 16-01-00409}%
{}\footnote{Steklov Mathematical Institute RAS,8 Gubkin St., Moscow 119991 Russia, and Department of Mathematics, National Research University the Higher School of Economics, Usacheva Str., 6 Moscow 119048 Russia, e-mail: \emph{kazarian@mccme.ru}}}

\title{Non-associative Hilbert scheme and Thom polynomials}
\date{}
\usepackage{amssymb,amsmath}
\usepackage{theorem,graphics}
\usepackage{longtable}
\usepackage
{hyperref}
\usepackage[all]{xy}

\newcounter{ssect}[section]
\renewcommand{\thessect}{\thesection.\arabic{ssect}}
\newcommand{\ex}[1][?]{\refstepcounter{ssect}%
\def\thnum{\ifx?#1\else\label{#1}\fi{ \thessect}}\exold}

\def\mdskip{\vskip-\lastskip\vskip\medskipamount}
\long\def\exold#1.#2\par{\ifvmode\mdskip\fi{\bf#1\thnum.}#2\par\mdskip}
\long\def\th#1.#2\par{\ex#1.{{\em#2}}\par}
\def\QED{\ifmmode\eqno\square\else
{\parfillskip0pt\hfil$\square$\par}\mdskip\fi}

\def\C{{\mathbb C}}

\def\Z{{\mathbb Z}}
\def\a{{\alpha}}

\def\cO{{\cal O}}

\def\deg{\mathop{\rm deg}}
\def\Hom{{\rm Hom}}

\def\codim{\mathop{\rm codim}\nolimits}

\def\pt{{\rm pt}}

\let\l\lambda

\def\rk{{\rm rk}}
\def\*{{\otimes}}

\def\Sym{\mathop{\rm Sym}\nolimits}

\let\o\omega
\def\res{\mathop{\rm res}\limits}

\def\@{{\odot}}

\def\Alg{{\rm Alg}}
\def\fm{{\mathfrak m}}
\def\GL{{\rm GL}}
\def\GLN{{\rm GL}(\ud)}
\let\g\gamma
\let\<\langle
\let\>\rangle

{\theorembodyfont{\rmfamily}    

}

\let\S\Sigma
\newcommand{\too}[1][\relax]{\stackrel{#1}\longrightarrow}
\newcommand{\oot}[1][\relax]{\stackrel{#1}\longleftarrow}
\renewcommand{\to}[1][\relax]{\stackrel{#1}\rightarrow}

\def\ud{{\underline d}}
\def\Hilb{{\rm Hilb}}
\def\loc{{\rm loc}}
\def\tK{{\widetilde K}}
\def\Tp{{\rm Tp}}
\let\b\beta
\let\k\varkappa
\def\cK{{\mathcal K}}
\def\Id{{\rm Id}}

\begin{document}

\maketitle

\begin{abstract}
Thom polynomial describes the cohomology class Poincar\'e dual to the locus of particular singularity of a generic holomorphic map. In this paper we derive a closed formula for the generating function of its coefficients. The method is based on a new construction of the embedding space of punctual Hibert scheme that we call the non-associative Hilbert scheme. The efficiency of the method is demonstrated on explicit computation of a number of Thom polynomials, including those associated with singularities of Thom-Boardman types~$\S^{i,j}$, $\S^{2,1,1}$, $\S^{2,2,1}$, and~$\S^{2,2,2}$.
\end{abstract}

\section{Introduction}

Thom polynomial is associated with a local singularity type of
holomorphic mappings. If $f:X\to Y$ is a generic holomorphic map
between two nonsingular varieties then the Thom polynomial
evaluated on the Cherm classes of these varieties expresses the
cohomology class Poincar\'e dual to the locus in~$X$ of a given
singularity type for~$f$. The theory of Thom polynomials is
deserved as a possible approach to some problems of enumerative
projective algebraic geometry~\cite{DS6,Kaz03}. Some interesting applications
were discovered recently in another domains of mathematics such
that representation theory and the study of moduli spaces~\cite{FNR08,FNR08a,KL1,KL2}.

The computation of Thom polynomials for particular singularity
classes was initiated in 60-ies by R.Thom, and many researches
contributed to this activity, see~\cite{Thom,Por,Ro,D,Gaf81,Rim01}. It was observed quite
early~\cite{Ro,D} that in many cases, for example, for contact
classification of singularities, the Thom polynomial is expressed
in terms of the relative Chern classes $c_i=c_i(f^*TY-TX)$ of the
map. Moreover, it is independent of the dimensions of the
manifolds~$X$ and~$Y$ provided that the relative dimension
$\ell=\dim Y-\dim X$ is fixed. If by the singularity type we mean
the isomorphism class of the local algebra then it is possible to
relate singularity classes for different values of
the relative dimension~$\ell$. The corresponding Thom polynomials
differ dramatically, and even its degree depends on~$\ell$.

However, it was noticed by Feh\'er and Rim\'anyi in~\cite{FRthser}, that
 it is possible to
represent the Thom polynomial in the form of an infinite series
(aka \emph{Thom series})
\begin{equation*} \Tp_\eta=
   \sum_{i_1,\dots,i_{\mu}}
   \psi_{i_1,\dots,i_{\mu}}c_{\ell+1+i_1}\dots
   c_{\ell+1+i_{\mu}}
\end{equation*}
whose integer coefficients $\psi_{i_1,\dots,i_{\mu}}$ are
independent of~$\ell$. Here the singularity type~$\eta$ is assumed
to be determined by the isomorphism class of the corresponding
local algebra and $\mu$ is the dimension of the local algebra
diminished by~$1$. First computations of Thom series can be found
in~\cite{BFR03,FRthser,FRloc,Pr05,Pra07}.

A new achievement in the development of the Thom polynomial theory was
attained in~\cite{BS}. In that paper the authors applied the
Atyah-Bott's localization principle to the computation of the Thom
polynomials for singularities of Thom-Boardman type
$\S^{1_d}=\S^{1,\dots,1}$. The final formula obtained in that
paper can be reformulated by saying that the coefficients of the
Thom series are determined by a rational generation function whose
terms are identified explicitly for $d\le 6$ and up to certain
undetermined coefficients for bigger~$d$. The computations of this
paper didn't use the existence theorem for Thom polynomials and
the very fact that the result is expressed in terms of the
relative Chern classes is obtained by a lengthy direct
computations with unexpected magic cancelations.

The localization methods was extended in~\cite{FRloc} to the study
of Thom polynomials for other kinds of singularities. It is
established, in particular, that the whole Thom series can be
recovered from a finite number of its initial terms. In other
words, the knowledge of Thom polynomial for some $\ell$ big enough
(in fact, for $\ell=\mu$) is sufficient to recover this polynomial
for all~$\ell$.

In the present paper we summarize the mentioned above results on
the Thom series by presenting an explicit rational form of the
generating function for its coefficients (Theorem~\ref{th1}). Our
computation is direct and does not use an \emph{apriori} existence
theorem for Thom polynomials. This allows us to use minor
assumptions on the chosen singularity type: it could be given by
any algebraic family of finite dimensional local algebras. For the
particular case of the singularity types $\S^{1_d}$ the formula of
Theorem~\ref{th1} specializes to that from~\cite{BS}. However, our
proof is essentially different.

The computation uses an old construction of partial resolution of
the singularity locus by means of the punctual Hilbert scheme. The main
novelty is a new geometric construction for the smooth embedding
space of the Hilbert scheme that we call a nonassociative Hilbert
scheme. The embedding of the Hilbert scheme to its ambient space
can be identified as the solution locus of the associativity
equation. The computational part of the proof of Theorem~\ref{th1}
is the application of the Gysin homomorphism for the constructed
partial resolution. At this point we use the formalism developed
in~\cite{Kaz11}. It allows one to keep the answer expressed in
terms of the relative Chern classes at all intermediate steps.
However, whenever the geometric construction is presented, the
computation of the Gysin homomorphism is technical, and any other
known ways of its computation, including the Schubert calculus or
the localization formulas, would lead to the same answer.

All terms of the rational function of Theorem~\ref{th1} are
identified explicitly except certain polynomial $P_\eta$ whose
computation is not covered by the assertion of the theorem. We
call this polynomial the local characteristic invariant of the
singularity~$\eta$. It is defined as the equivariant Poincar\'e dual
class of certain explicitly given subvariety in a vector space.
For the singularity determined by fixing the dimension vector of
local algebra, it is given by the associativity equation.

At the moment, we have no regular efficient way to compute $P_\eta$ in general, but first examples of its computation are also presented in this paper. It leads to a considerable list of closed formulas for Thom polynomials of different singularities. It includes also some cases for which the Thom polynomial was unknown before, in particular, those of Thom-Boardman types~$\S^{2,1,1}$, $\S^{2,2,1}$, and~$\S^{2,2,2}$.

In the present paper, we use the topological language of complex geometry and cohomology. However, all results are also valid in the algebraic context of arbitrary algebraically closed ground field of zero characteristic and the Chow rings instead of cohomology. When necessary, we supply technical details allowing one to adopt our topological arguments to the algebraic setting.

The paper is organized as follows. In Section~\ref{secTh} we
review the notions of Thom polynomials and equivariant Poincar\'e
duals, compare topological and algebraic aspects of their
definitions and revise the proof of the existence theorem for Thom
polynomials. The main Theorem~\ref{th1} expressing the rational generating function for the Thom series is formulated in section~\ref{sec3}.
The definition of the polynomial~$P_\eta$ used in the formulation of Theorem~\ref{th1} is introduced in Sect.~\ref{sectAss}. In Section~\ref{sec5} we
recall the resolution method for the computation of Thom
polynomials. The main geometric construction of the paper---the
nonassotiative Hilbert scheme---is introduced in Sect.~\ref{sec6}.
The computation of Thom polynomial is completed in
Section~\ref{sec7} by applying the Gysin homomorphism. The results on explicit computations based on Theorem~\ref{th1} for particular singularity classes are put together in Sect.~\ref{sec8}. In the final Section~\ref{sec9} we formulate some open problems.

\textbf{Acknowledgments.} My sincere gratitude is to L. Feh\'er, R. Rim\'anyi, and A. Szenes. The paper appeared as the result of long clarifying and stimulating discussions with them.

\section{The notion of Thom polynomial}\label{secTh}

\subsection{Equivariant Poincar\'e duals and multidegrees}

Let a complex algebraic Lie group $G$ act algebraically on a vector
space~$V$. Then for any algebraic $G$-invariant subvariety
$\eta\subset V$ of pure codimension $k$ there is a well defined its
equivariant Poincar\'e dual cohomology class
$$[\eta]\in H^{2c}_{G}(V)=H^{2k}_{G}(\pt)=H^{2c}(BG).$$

If $G=T^n=(\C^*)^n$ is the (complex) torus, then
$H^*(BG)=\Z[t_1,\dots,t_n]$ with the (real) grading of the class
$t_i$ equal to~$2$ for all $i=1,\dots,n$. The polynomial
representing~$[\eta]$ in this case is also known in algebraic
geometry as the \emph{Josef polynomial} or \emph{multidegree}. If $f:V\to\C$ is a $T^n$-homogeneous function satisfying
$$f(\l v)=\l_1^{a_1}\dots\l_n^{a_n}\,f(v),
 \qquad\l=(\l_1,\dots,\l_n)\in T^n,$$
then the multidegree of its zero hypersurface is equal to
$$[\{f=0\}]=a_1t_1+\dots+a_nt_n.$$

The multidegree of a (reduced) complete intersection given as the zero set of a collection of homogeneous functions is equal to the product of the corresponding linear factors. This can be applied, in particular, to coordinate subspaces. The multidegree of any $T^n$-invariant subvariety can be computed by deforming it in a flat family to a union of coordinate subspaces (probably, with multiplicities). Then its multidegree is equal to the corresponding linear combination of multidegrees of coordinate subspaces.

In the case when $G=\GL(n)$ we have $BG=G_{n,\infty}$, the
infinite Grassmann manifold, and $H^*(BG)$ is the polynomial ring
generated by the Chern classes $c_i$ of the tautological vector
bundle over the Grassmannian. Therefore, the equivariant
Poincar\'e dual class $[\eta]$ in this case is a polynomial with
integer coefficients in $c_1,\dots,c_n$ of weighted degree~$k$
where we assume $\deg(c_i)=i$. Restricting the action of $\GL(n)$
to the maximal torus $T^n\subset\GL(n)$ we can express $[\eta]$ as a
polynomial in the Chern roots $t_1,\dots,t_n$. This polynomial is
symmetric and~$c_i$ corresponds to the~$i$th elementary symmetric
function in the variables~$t_1,\dots,t_n$.

%
%
%

\subsection{Thom polynomials}

A \emph{Thom polynomial} $\Tp_\eta$ is a special case of the
equivariant Poincar\'e dual class when $V=J^K_0(\C^m,\C^n)$ is the
space of $K$-jets of map germs from $\C^m$ to $\C^n$ at the
origin, $G$ is the group of $K$-jets of changes of coordinates in the source and the target, respectively, and $\eta$ is an arbitrary singularity type. By a `singularity type' we mean here any $G$-invariant reduced irreducible algebraic subvariety in~$J^K_0(\C^m,\C^n)$. The group~$G$ is homotopy equivalent to its subgroup $\GL(m)\times \GL(n)$ consisting of linear changes. Therefore, we have, by definition,
$$\Tp_\eta=[\eta]\in H_G^*(V)=H^*_{\GL(m)\times\GL(n)}(\pt)
 =H^*(B\GL(m)\times B\GL(n))=\Z[c_1',\dots,c_m',c_1'',\dots,c_n''].$$

We see, that from the very definition, $\Tp_\eta$ is a polynomial in two groups of variables, $c_i'$ and $c_j''$ corresponding to the source and the target, respectively.

The definition depends on a choice of the order $K$ of jets. We say that the singularity~$\eta\subset J^K(\C^m,\C^n)$ is \emph{$k$-determined} if it is of the form $\eta=p^{-1}(\eta_k)$ where $\eta_k\subset J^k_0(\C^m,\C^n)$ and $p:J^K(\C^m,\C^n)\to J^k(\C^m,\C^n)$ is the natural projection. The Thom polynomial of a $k$-determined singularity is actually independent of a choice of~$K\ge k$.

The importance of the notion of Thom polynomial
shows up in the following theorem~\cite{Thom,HK,Kaz00}. Let $X$ and $Y$ be manifolds of dimensions $\dim X=m$, $\dim Y=n$. The jet bundle $J^K(X,Y)$ is formed by all $K$-jets of maps from $X$ to $Y$ at all points $x\in X$ and $y\in Y$. It is the total space of a bundle over $X\times Y$ with the fiber isomorphic to $J^K_0(\C^m,\C^n)$. It follows that the projection $J^K(X,Y)\to X\times Y$ is a homotopy equivalence and the induced homomorphism in cohomology is an isomorphism. Denote by $\eta(X\times Y)\subset J^K(X,Y)$ the subvariety formed by jets of maps with prescribed singularity~$\eta$.

\th[Thom] Theorem. The cohomology class
$$[\eta(X\times Y)]\in H^*(J^K(X,Y))=H^*(X\times Y)$$
represented by the locus of the singularity type  $\eta$ is equal to the Thom polynomial $\Tp_\eta$ evaluated on the Chern classes of the tangent bundles $c_i(TX)$ and $c_j(TY)$ of the factors in $X\times Y$.

\th Corollary. Let $f:X\to Y$ be a holomorphic map, $j^Kf:X\to J^K(X,Y)$ be its jet extension, and $\eta(f)=(j^Kf)^{-1}\eta(X\times Y)\subset X$, the locus of the singularity~$\eta$ for the map~$f$. Then, if $\eta(f)$ has expected dimension and reduced then its fundamental class in~$X$ represents the Thom polynomial $\Tp_\eta$ evaluated on the classes $c_i(TX)$ and $f^*c_j(TY)$. For an arbitrary map, $\Tp_\eta$ is represented by an algebraic cycle supported on~$\eta(f)$.

\medskip
The notions of Chern classes, equivariant
Poincar\'e duals, and Thom polynomials make sense in the algebraic
geometry context with Chow groups instead of cohomology.
The assertions of the theorem and its corollary remain true in the algebraic setting. Since the proof of this generalization has never been published we reproduce it here. The topological counterpart of the proof
uses the notions of classifying spaces and classifying maps for $G$-bundles. We will see that these notions can be adopted to the algebraic context, see also~\cite{Tot}.

Let us first review the definition of the Chern classes of a given vector bundle~$E$ over a nonsingular algebraic base~$X$. In topology, the Chern classes are certain distinguished elements in the cohomology of the Grassmannian that serves as the classifying space for vector bundles. Since the Chow ring of the Grassmannian is isomorphic to the cohomology, in order to define Chern classes, it is sufficient to construct a classifying map from~$X$ to the Grassmannian that induces the given vector bundle~$E$. Here is the construction of such classifying map adopted to the algebraic setting. Choose some integer $N\gg 0$ an define $X_1$ to be the total space of the bundle $\Hom(E,\C^N)$ over~$X$. Then obviously $H^*(X)\simeq H^*(X_1)$ and the same holds for Chow groups. The point of $X_1$ is a linear map from a fiber of $E$ to $\C^N$. Let $X_2\subset X_1$ be the open subset consisting of injective linear maps. The complement $X_1\setminus X_2$ has a very big codimension, growing together with~$N$, therefore, the inclusion $X_2\to X_1$ induces an isomorphism in cohomology for any fixed grading, $H^k(X_2)\simeq H^k(X_1)$, if $N$ is big enough, and a similar isomorphism holds in Chow groups.

The image of an injective linear map can be regarded as a point of the Grassmannian. This gives the required classifying map $\k:X_2\to G_{n,N}$, $n=\rk E$. Then we define the Chern classes~$c_j(E)$ via the induced homomorphism
$$\k^*:H^{2j}(G_{n,N})\to H^{2j}(X_2)\simeq
  H^{2j}(X_1)\simeq H^{2j}(X).$$
In the algebraic setting the Chern classes defined in this way (for a smooth base~$X$) are classes of rational equivalence of $j$-codimensional subvarieties.

Now we present a similar construction for jet bundles. For given $m$, $n$, and $K$ chose some $N\gg0$.

\ex Definition. \emph{A $K$-jet of an $m$-dimensional submanifold} in $\C^N$ is the $K$-jet of a map germ $(\C^m,0)\to(\C^N,0)$ with the derivative of full rank at the origin and considered up to a change of coordinates in the source~$\C^m$. The variety of all $k$-jets of submanifolds is denoted by $\widetilde G_{m,N}$.

Assignment of the tangent plane to a submanifold defines a natural mapping~$\widetilde G_{m,N}\to G_{m,N}$. This mapping can be represented as a tower of fibrations with affine fibers. This determines the structure of a nonsingular algebraic variety on~$\widetilde G_{m,N}$ and also shows an isomorphism between Chow groups and cohomology of the varieties $\widetilde G_{m,N}$ and $G_{m,N}$.

\ex Definition. \emph{A $K$-jet of an $n$-dimensional `quotient' manifold} of $\C^N$ is the $K$-jet of a map germ $(\C^N,0)\to(\C^n,0)$ with the derivative of full rank at the origin and considered up to a change of coordinates in the target~$\C^n$. The variety of all $K$-jets of quotient manifolds is denoted by $\widehat G_{N-n,N}$.

Similarly to the case of submanifolds, we have a natural projection $\widehat G_{N-n,N}\to G_{N-n,N}$ that induces an isomorphism both in cohomology and in the Chow groups.

\ex Definition. The \emph{classifying space of singularities} is defined to be the product space
$$B=B_{m,n,N}^K=\widetilde G_{m,N}\times \widehat G_{N-n,N}.$$

Again, the natural projection $B\to G_{m,N}\times G_{N-n,N}$ induces an isomorphism both in cohomology and in the Chow groups.

A point $(f,g)\in B$ determines a $K$-jet of a map germ $g\circ f:\C^m\to \C^n$ defined up to a change of coordinates in the source and the target. Consider the subvariety $\eta(B)\subset B$ formed by points with the associated map of prescribed singularity type~$\eta$. The following definition is equivalent to the preceding one.

\ex Definition. The \emph{Thom polynomial} of the singularity $\eta$ is the class represented by the subvariety $\eta(B)\subset B$ in either cohomology or the Chow ring of $B\sim G_{m,N}\times G_{N-n,N}$ and expressed in terms of the multiplicative generators of this ring for which we take traditionally the Chern classes of the tautological rank~$m$ subbundle and the rank~$n$ quotient vector bundles over the product of the Grassmannians.

It is easy to see that the Thom polynomial is independent of~$N$ provided that~$N$ is big enough.

Let now $X$ and $Y$ be manifolds of dimensions $m$ and $n$, respectively. Set
$$X_1=J^K(X,\C^N),\qquad Y_1=J^K(\C^N,Y),$$
and consider open subsets $X_2\subset X_1$ and $Y_2\subset Y_1$ formed by jets of maps with the derivative of full rank.

On one hand, the composition of maps defines a natural surjective mapping
$$\pi:X_2\times Y_2\to J^K(X,Y).$$
This mapping induces an isomorphism in the cohomology or the Chow groups of fixed grading provided that $N$ is big enough. Therefore, we can replace the class of the subvariety $\eta(X\times Y)\subset J^K(X,Y)$ participating in the statement of Theorem~\ref{Thom} by the class of its preimage $\pi^{-1}(\eta(X\times Y))$ in $X_2\times Y_2$.

On the other hand, the classifying `tautological' map $\k:X_2\times Y_2\to B$ is well defined. This map respects singularity types, $\k^{-1}(\eta(B))=\pi^{-1}(\eta(X\times Y))$ and the pullbacks of the basic Chern classes under $\k^*$ are the Chern classes of tangent bundles of~$X$ and~$Y$, respectively, see the diagram.
$$\xymatrix{
H^k(X\times Y)\too[\simeq]H^k(J^K(X,Y))\too[\simeq]
   H^k(X_2\times Y_2)&&{H^k(B)
 \makebox[0pt][l]{${}
   \oot[\simeq]H^k(G_{m,N}\times G_{N-n,N})$}}
   \ar[ll]_-{\k^*}\\
   [\eta(X\times Y)]\ar@{}[u]|-{\rotatebox{90}{$\in$}}
   &&[\eta(B)]
 \makebox[0pt][l]{${}
   =\Tp_\eta$}
   \ar@{|->}[ll]\ar@{}[u]|-{\rotatebox{90}{$\in$}}
}\makebox[11em]{}$$

This proves both Theorem~\ref{Thom} and its algebraic geometry counterpart.\QED

\subsection{Thom polynomial of a contact singularity}

A contact singularity is a particular case of a singularity type and a general definition of Thom polynomial is applicable in this case as well. There are, however, some details that we are going to clarify here.
The \emph{local algebra} of a map germ $f:(\C^m,0)\to(\C^{n},0)$ is
defined to be the quotient of the ring of functions
$\C[[x_1,\dots,x_m]]$ on the source modulo the ideal generated by
the components $f_1,\dots, f_{n}$ of the mapping,
$$Q_f=\frac{\C[[x_1,\dots,x_m]]}{(f_1,\dots,f_{n})}.$$

Two map germs are called \emph{contact equivalent}, or $\cK$-\emph{equivalent}, if their local algebras are isomorphic.

We will assume that $\dim Q_f<\infty$. If $n\ge m$, then the map germs with infinite dimensional local algebras form a subset of infinite codimension in the space of map germs. We do not consider such maps in this paper.

In a naive approach, a contact singularity type is the collection of map germs with a fixed isomorphism class of local algebras. This possible definition has a disadvantage. Consider, for example, the Porteous-Thom singularity $\S^r$ defined by the requirement that the derivative of the map has an $r$-dimensional kernel. It is natural to associate the local algebra~$Q_{\S^r}=\C[x_1,\dots,x_r]/\fm^2$ to this singularity, where~$\fm$ is the maximal ideal. The local algebra of a generic map germ with the  singularity~$\S^r$ is indeed isomorphic to $Q_{\S^r}$, if the relative dimension $\ell=n-m$ of the map satisfies $\ell\ge \frac{r(r+1)}{2}$. If $\ell<\frac{r(r+1)}{2}$ then this local algebra is never realized since $m+\frac{r(r+1)}{2}$ is the smallest number of generators of the ideal in $\C[[x_1,\dots,x_m]]$ whose quotient algebra is isomorphic to~$Q_{\S^r}$. However, it makes sense to consider the singularity $\S^r$ and its Thom polynomial for any~$\ell\ge-r$, even including some negative values of the relative dimension.

Let $Q$ be a finite dimensional local algebra and $f:(\C^m,0)\to(\C^n,0)$ be a map germ.

\ex[Ksing] Definition. {We say that \emph{$f$ has singularity type $Q$} at the origin if there exist at most finitely many ideals~$I$ in $\C[[x_1,\dots,x_m]]$ satisfying two requirements: 1) the quotient algebra $\C[[x_1,\dots,x_m]]/I$ is isomorphic to~$Q$, and 2) the ideal generated by the components $f_j$ of the germ~$f$ is contained in~$I$. The number of ideals satisfying these requirements is called the \emph{multiplicity} of the singularity~$Q$ at the germ~$f$.

By a \emph{contact singularity type} we mean now an isomorphism class of algebras or an (irreducible) algebraic family of local algebras of fixed dimension.

The \emph{Thom polynomial of a contact singularity} is the one associated with the closure of map germs of prescribed contact singularity type, multiplied by the multiplicity of its generic representative.}

If $\ell$ is big enough (if $m+\ell$ is bigger than the minimal number of generators of an ideal $I\subset \C[x_1,\dots,x_m]$ with the quotient algebra isomorphic to~$Q$) then for a generic germ with all components in $I$ the components do generate the ideal~$I$. In this case a generic germ with the singularity $Q$ has local algebra isomorphic to~$Q$ and the multiplicity is equal to~$1$. In other words, a special account is required for the case of `small~$\ell$' only.

It is useful to note that the contact singularity class associated with any finite dimensional local algebra~$Q$ is finitely determined, namely, it is at least $\mu$-determined where $\mu=\dim Q-1$. This is a straightforward consequence of the fact that any ideal of finite codimension $\mu+1$ in the local coordinate ring $\C[[x_1,\dots,x_m]]$ contains the $(\mu+1)$st power of the maximal ideal. It follows that the condition on a given map germ to have a contact singularity of type~$Q$ is uniquely determined by the $\mu$-jet of the germ and that the corresponding Thom polynomial is independent of the order~$K$ of jet space used in its definition provided that $K\ge\mu$.

\qquad
The Thom polynomial of a contact singularity is expressed actually as polynomial in the \emph{quotient variables} $c_i$ defined by
$$1+c_1+c_2+\dots=\frac{1+c_1''+\dots+c_n''}
   {1+c_1'+\dots+c_m'}.$$
For a given holomorphic map $f:X\to Y$ the quotient variables specialize to the relative Chern classes $c_i=c_i(f)=c_i(f^*TY-TX)$. Moreover, the Thom polynomial is independent of the dimensions $m$ and $n$ of the source and the target of the map provided that the relative dimension $\ell=n-m$ is fixed. This fact has a simple a priori explanation, see~\cite{Ro,D,Kaz00}. On the other hand, it is accounted automatically in the main formula of the next section for Thom polynomials of contact singularities.

\section{Main theorem}\label{sec3}

It will be more convenient for us to encode contact singularities by isomorphism classes of finite dimensional \emph{nilpotent algebras}. The nilpotent algebra~$N_f$ of a map germ $f:(\C^m,0)\to(\C^n,0)$ is the maximal ideal in the corresponding local algebra,
$$N_f=\frac{\fm_{x_1,\dots,x_m}}{(f_1,\dots,f_n)},$$
where $\fm_{x_1,\dots,x_m}$ is the maximal ideal in the local coordinate ring $\C[[x_1,\dots,x_n]]$. The nilpotent algebra carries the same information as the local one; the algebra $Q_f$ can be reconstructed from $N_f$ by just adding unity. In particular, $\dim Q_f=\dim N_f+1$.

Let $N$ be a finite dimensional nilpotent algebra. Choose a filtration on it, that is, a finite decreasing sequence of subalgebras
$$N=N_1\supset N_2\supset\dots\supset N_{r+1}=0$$
satisfying $N_i\cdot N_j\subset N_{i+j}$.
The dimension of the subsequent quotients of this
filtration is denoted by $d_k=\dim(N_k/N_{k+1})$. The sequence $\ud=(d_1,\dots,d_r)$ is
called the \emph{dimension vector}. Set also
$$\mu=d_1+\dots+d_r=\dim N.$$

A filtration on the nilpotent algebra always exists but it is usually not unique. A formula for the Thom polynomial presented below depends on a choice of the filtration (although the Thom polynomial itself is obviously independent of this choice). A particular filtration is given by the powers of the maximal ideal, $N_k=N^k$. In this case all components $d_i$ of the dimension vector are strictly positive. But in general, we allow for some components of the dimension vector to be vanishing. (I am grateful to R.Rimanyi who payed my attention to the possibility of using filtrations different from those formed by powers of the nilpotent algebra).

For an integer $i$ such that $d_1+\dots+d_{k-1}< i\le
d_1+\dots+d_k$, we define its \emph{weight} to be $w(i)=k$ and the
\emph{exponent} to be $e(i)=d_1+\dots+d_{k}-i$. The weights and
the exponents of the numbers $1,2,\dots,\mu$ form two sequences
\begin{equation}
\begin{gathered}
\underline
 w=(\underbrace{1,\dots,1}_{d_1},\underbrace{2,\dots,2}_{d_2},\dots,
 \underbrace{r,\dots,r}_{d_r}),\\
 \underline e=(\underbrace{d_1-1,\dots,1,0}_{d_1},
  \underbrace{d_2-1,\dots  ,1,0}_{d_2},\dots,
  \underbrace{d_{r}-1,\dots,1,0}_{d_r}).
\end{gathered}
 \label{we}\end{equation}
Introduce the following rational function determined uniquely by the dimension vector
\begin{equation}
K_{\ud}(t_1,\dots,t_\mu)=
   \frac{\displaystyle\prod_{i=1}^\mu
   t_i^{e(i)}\;\prod_{1\le i<j\le\mu}\!\!\!(t_j-t_i)}
   {\displaystyle\prod_{\substack{1\le i\le j<k\le\mu\\w(i)+w(j)\le
   w(k)}}\!\!\!\!\!\!(t_k-t_i-t_j)}.
 \label{eqK}\end{equation}

Consider any contact singularity type $\eta$ given as a unique isomorphism class of nilpotent algebras or given by an irreducible algebraic family of nilpotent algebras with fixed dimension vector~$\ud$.

\th[th1] Theorem. The Thom polynomial of the singularity $\eta$
with the given dimension vector $\ud=(d_1,\dots,d_r)$ can be
represented in the form of a series
\begin{equation}\label{Thser}
\Tp_\eta=
   \sum_{i_1,\dots,i_{\mu}}
   \psi_{i_1,\dots,i_{\mu}}c_{\ell+1+i_1}\dots
   c_{\ell+1+i_{\mu}}
\end{equation}
in the quotient variables {\rm(}the relative Chern classes{\rm)} $c_1,c_2,\dots$ whose coefficients $\psi_{i_1,\dots,i_{\mu}}$ are
independent of $\ell$ and are determined by the following rational
generating function
$$\sum_{i_1,\dots,i_{\mu}}
   \psi_{i_1,\dots,i_{\mu}}t_1^{i_1}\dots t_{\mu}^{i_{\mu}}=
   K_{\ud}(t_1,\dots,t_\mu)\cdot P_\eta(t_1,\dots,t_\mu),$$
where~$P_\eta$ is a polynomial defined below in
Sect.~{\rm\ref{sectAss}}.

By the power series expansion of a rational function we mean its expansion in the domain $t_1\ll t_2\ll\dots\ll t_\mu$. It means the following. First, we expand the function as a Laurent series in $t_\mu^{-1}$ whose coefficients are rational functions in $t_1,\dots,t_{\mu-1}$. Then we expand each of the coefficients as Laurent series in $t_{\mu-1}^{-1}$, and so on.

Equivalently, the relation of Theorem~\ref{th1} can be represented formally in the form of iterated residue~\cite{BS}
$$\res_{t_1=\infty}\dots\res_{t_\mu=\infty}
K_{\ud}(t)\cdot P_\eta(t) \prod_{i=1}^\mu\Bigl(-C(1/t_i)\,t_i^{\ell}dt_i\Bigr),\quad
C(1/t)=\sum_{k=0}^\infty c_k t^{-k}.$$

The actual definition of the polynomial $P_\eta$ is discussed in the next section. Even without knowing explicitly its coefficients we can formulate some general conclusions resulting from its existence. The fact that the Thom polynomial can be represented in a form of a series~\eqref{Thser} was observed first by L.Feh\'er and R.Rim\'anyi  in~\cite{FRthser}. Remark that the whole series is usually infinite but it has only finitely many nonzero terms for each particular value of~$\ell$. Every monomial in the Thom polynomial has homogeneous degree at most $\mu$, and the change of the relative dimension~$\ell$ is accounted by a shift of the indices in the Chern classes. It follows that the codimension of the singularity~$\eta$ satisfies
$$\codim\eta=(\ell+1)\,\mu+d$$
where $d$ is independent of~$\ell$. We can regard $d$ as the `virtual codimension' of the singularity~$\eta$ in the case of the relative dimension~$\ell=-1$.

It was proved in~\cite{FRloc} that the whole Thom series can be formally reconstructed from finitely many its terms. The statement of the theorem makes this reconstruction procedure more explicit: we see that the series is uniquely determined by finitely many coefficients of the polynomial $P_\eta$.

The language of generating series in presentation of characteristic classes was introduced in~\cite{Kaz11} to simplify some formal manipulations with them. For example, it allows one to pass easily from the basis of Chern monomials to the basis consisting of Schur functions. For any sequence of integers $\l=(\l_1,\dots,\l_\mu)$ we set
$$\Delta_{\l_1,\dots,\l_\mu}=\det|c_{\l_i-i+j}|_{1\le i,j\le \mu}.$$
If the sequence $\l$ is non-increasing, then this is the Schur function associated with the partition~$\l$. For any sequence, not necessarily non-increasing, this determinant is either zero or is equal to some Schur function, up to a sign. In terms of generating series passing from the basis of Chern monomials to the basis of Schur functions is accounted by a simple factor, see~\cite{Kaz11} for the proof.

\th[cor1] Corollary. The Thom polynomial admits a Schur function expansion
$$\Tp_\eta=\sum
   \sum_{i_1,\dots,i_{\mu}}
   \widetilde \psi_{i_1,\dots,i_{\mu}}
   \Delta_{\ell+1+i_1,\dots,\ell+1+i_{\mu}}
$$
with the following generating function for its coefficients
$$\sum_{i_1,\dots,i_{\mu}}
   \widetilde\psi_{i_1,\dots,i_{\mu}}t_1^{i_1}\dots t_{\mu}^{i_{\mu}}=
   \tK_{\ud}(t_1,\dots,t_\mu)\cdot P_\eta(t_1,\dots,t_\mu),$$
where $P$ is the same as in Theorem~\ref{th1} and
\begin{equation}
\tK_\ud=\frac{K_\ud}{\displaystyle\prod_{1\le i<j\le\mu}\!\!\!\bigl(1-\frac{t_i}{t_j}\bigr)}
    =   \frac{\displaystyle\prod_{i=1}^\mu
   t_i^{e(i)+i-1}}
   {\displaystyle\prod_{\substack{1\le i\le j<k\le\mu\\w(i)+w(j)\le
   w(k)}}\!\!\!\!\!\!(t_k-t_i-t_j)}.
 \label{eqKt}\end{equation}

%

\section{Definition of the polynomial~$P_\eta$ of Theorem~\ref{th1}}\label{sectAss}

The polynomial~$P_\eta$ participating in Theorem~\ref{th1} is the product of two factors,
$$P_\eta=P'_\eta\cdot P''_\eta.$$
We call these factors the local characteristic invariant of the singularity~$\eta$ and the normalizing factor, respectively. The definition of these factors is given below.

\subsection{Associativity equation}

Consider a flag of finite dimensional vector spaces
$$N=N_1\supset N_2\supset\dots\supset N_{r+1}=0,\qquad \dim N_{k}/N_{k+1}=d_k.$$
A \emph{filtered commutative algebra structure} on $N$ is a linear
mapping
$$\psi:\Sym^2N\to N$$
such that $\psi(N_k\*N_m)\subset N_{k+m}$. The filtered
commutative algebra structures on $N$ form a vector space that we
denote by $\Alg(\ud)=\Alg(d_1,\dots,d_r)$.

Set $\mu=\dim N=\sum_{i=1}^r d_i$ and choose a basis $e_1,\dots, e_\mu$ in $N$ such that the vectors $e_{d_1+\dots+d_{k-1}+1},\dots,e_{d_1+\dots+d_{k-1}+d_k}$ generate $N_k/N_{k+1}$. Then the multiplication law is given by
$$e_i\cdot e_j=\psi(e_i,e_j)=\sum_{k,w(k)\ge w(i)+w(j)} q_{i,j}^k e_k,$$
where the weight sequence $(w(1),\dots,w(\mu))$ is determined by the dimension vector $(d_1,\dots,d_r)$ as in~\eqref{we}. The coefficients
$$q_{i,j}^k=q_{j,i}^k,\quad w(k)\ge w(i)+w(j),$$
form the coordinate system in the vector space $\Alg(\ud)$.

The multiplication operation does not satisfy automatically the
associativity law. Denote by $\Hilb_{\ud}^\loc\subset\Alg(\ud)$ the
subvariety consisting of associative algebras. It is easy to see
that it is algebraic. Moreover, its ideal is generated by
quadratic equations of the form
\begin{equation}
 \sum_{m}q_{i,j}^m q_{m,k}^n
    -\sum_{m}q_{i,m}^n q_{j,k}^m=0
 \label{eqass}\end{equation}
for all quadruples $(i,j,k,n)$ with $w(i)+w(j)+w(k)\le w(n)$. This identity
expresses the equality of the $n$-components of the products
$(e_i\cdot e_j)\cdot e_k$ and $e_i\cdot (e_j\cdot e_k)$.

Furthermore, if we are given a singularity
type $\eta$, denote by $\Hilb_\eta^\loc\subset\Hilb_\ud^\loc\subset
\Alg(\ud)$ the closure of the locus formed by the algebras
isomorphic to one of the nilpotent algebras defining~$\eta$. Both
subvarieties are invariant with respect to the natural action of
the group $\GLN$ of linear transformations of~$N$ preserving the
flag~$N_\bullet$. This group is homotopy equivalent to the group
$\GL(d_1)\times\dots\times \GL(d_r)$ and its ring of
characteristic classes is the polynomial ring generated by the
Chern classes of the $\GL(d_i)$-invariant vector spaces $N_i/N_{i+1}$,
$i=1,\dots,r$ (this is true both in the topological and in the
algebraic context). The maximal torus $(\C^*)^\mu\subset \GLN$ acts on $\Alg(\ud)$ by rescaling of the basic vectors $e_i$ in $N$. The multidegree of the coordinate $q_{i,j}^k$ with respect to this torus action is equal to
 $$\deg q_{i,j}^k=t_k-t_i-t_j.$$
Remark that the denominator of the function~$K_\ud$~\eqref{eqK} is nothing but the Euler class of $\Alg(\ud)$, that is, the product of multidegrees of all its coordinates.

\ex[lochi] Definition. The \emph{local characteristic invariant}
$P'_\eta(t_1,\dots,t_\mu)$ of the singularity~$\eta$ is defined to be the multidegree of the subvariety
$\Hilb_\eta^\loc\subset \Alg(\ud)$, that is,
its equivariant Poincar\'e dual cohomology class with respect to the natural action
of the group $\GLN\sim\prod_{i=1}^r \GL(d_i)=\prod_{i=1}^r
\GL(N_i/N_{i+1})$ and expressed in terms of the Chern roots of
that action.

\ex[exA4] Example. {If $r\le 2$, then the associativity equation is empty
and we have $\Hilb^\loc_\eta=\Alg(\ud)$ i.e. $P'_\eta=[\Hilb^\loc_\eta]=1$
if the chosen singularity class consists of all algebras with the
given dimension vector. The same is true if $r=3$ and the
dimension vector is $\ud=(1,1,1)$.

In the case $r=4$ and for the dimension vector $\ud=(1,1,1,1)$ we
obtain a nontrivial associativity restriction
$$q_{1,1}^2\cdot q_{2,2}^4-q_{1,3}^4\cdot q_{1,2}^3=0.$$
Generic algebras from the hypersurface defined by this equation is isomorphic to $A_4=\frac{x\C[[x]]}{(x^5)}$.
The multidegree of this hypersurface is equal to
$$P'_{A_4}(t_1,\dots,t_4)=t_{4}-t_{2}-2t_{1}.$$

In more complicated cases the subvariety $\Hilb_\eta^\loc\subset \Alg(\ud)$ is not a complete intersection and the polynomial $P'_\eta$ may contain
nonlinear factors. The way of its computation depends on a
particular choice of~$\eta$.}

It would be more convenient for us to regard below the algebra
structure on $N$ as the \emph{coalgebra} structure on the dual
space $N^\vee$ defined by the adjoint filtered morphism
$\psi^\vee:N^\vee\to \Sym^2 N^\vee$. The space $N^\vee$ is
equipped with the natural filtration
$$0=D_0\subset D_1\subset D_2\subset\dots\subset D_r=N^\vee,
 \qquad D_k={\rm ann} (N_{k+1})=(N/N_{k+1})^\vee,$$
and the comultiplication satisfies $\psi^\vee(D_k)\subset S_k$
where $S_k\subset\Sym^2N^\vee$ is spanned by the tensors $a\*b$,
$a\in D_i$, $b\in D_j$, with $i+j\le k$.

The subvariety $\Hilb_\ud^\loc\subset \Alg(\ud)$ is determined
by the associativity equation for this coalgebra and the
subvariety $\Hilb_\eta^\loc\subset  \Alg(\ud)$ is identified by
the condition that this coalgebra is isomorphic to the coalgebra
induced on the dual space $N_f^\vee$ to the nilpotent algebra of
one of the map germs~$f$ representing a given singularity type.

\subsection{Normalizing factor}

\ex[fnat] Definition. A filtration on a given nilpotent algebra~$N$ is called \emph{natural} if any automorphism of the algebra preserves also the filtration.

For example, the filtration formed by the powers of the algebra is natural. There could be, however, different filtrations, both natural and not. For a natural filtration we set $P''_\eta=1$.

If the filtration is not natural then there should exist a flag on $N$ forming a filtration different from the original one such that the algebra~$N$ with this filtration is isomorphic to the original one as a filtered nilpotent algebra. Denote by $F_\eta$ the closure of all flags with this property in the flag variety in~$N$.

\ex[norm] Definition. If the singularity~$\eta$ is defined by a unique isomorphism class of nilpotent algebras, we define $P''_\eta(t_1,\dots,t_\mu)$ as an arbitrary polynomial combination of Chern classes of the tautological vector bundles over the flag variety satisfying
$$\int_{F_\eta}P''_\eta(t_1,\dots,t_\mu)=1$$
and expressed in terms of the Chern roots of these bundles. If $\eta$ is given by a family of algebras, we apply the same definition to a generic member of the family.

By Chern roots we mean the formal substitution
$$c(N_1/N_2)=\prod_{i=1}^{d_1}(1+t_i),\quad\dots,\quad
  c(N_r/N_{r+1})=\prod_{i=d_1+\dots+d_{r-1}+1}^{d_1+\dots+d_r}(1+t_i).$$

A polynomial $P''_\eta$ satisfying this property does exist since any nonempty algebraic subvariety in the flag variety has a nonzero homology fundamental class. Moreover, there is a big freedom in a choice of this polynomial: possible choices of $P''_\eta$ form an affine hyperplane in the space of monomials of a fixed degree. This freedom can be used, for example, to achieve a cancelation with some factors in the denominator in the generating function of Theorem~\ref{th1}.

\section{Hilbert scheme and partial resolution of the singularity locus}\label{sec5}

Let $E$ and $F$ be finite dimensional vector spaces. The given singularity type $\eta$ constitutes a locus in the jet space of map germs between~$E$
and~$F$:
$$\eta\subset J^K_0(E,F)\simeq\Hom\Bigl(
 \bigoplus_{i=1}^K \Sym^{i}E,F\Bigr)=(\fm/\fm^{K+1})\otimes F$$
where $\fm=\fm_E$ is the maximal ideal in the ring of function germs on~$E$ at the origin. We need to compute the cohomology class $\Tp_\eta=[\eta]$ represented by this locus in the $(\GL(E)\times \GL(F))$-equivariant cohomology of the jet space. A step in the direction of this computation can be done using the following partial resolution of~$\eta$ the idea of which goes back to at least J.~Damon~\cite{D}.

Denote by $\Hilb_{\ud}(E)$ the variety parametrizing
flags of ideals
$$\fm=I_1\supset I_2\supset \dots\supset I_{r+1}$$
such that $\dim I_{k}/I_{k+1}=d_k$ for $k=1,\dots,r$ and $I_i\cdot
I_j\subset I_{i+j}$ for all~$i$ and~$j$ with $i+j\le r+1$. We
always have $\fm^{r+1}\subset I_{r+1}$, therefore, the variety
$\Hilb_{\ud}(E)$ can be embedded to an appropriate variety of flags in the finite dimensional vector space $\fm/\fm^{K+1}$ for any $K\ge r+1$. The variety $\Hilb_{\ud}(E)$ is equipped with the canonical subbundle $I=I_{r+1}/\fm^{K+1}$ of the trivial vector bundle $(\fm/\fm^{K+1})\times \Hilb_{\ud}(E)\to \Hilb_{\ud}(E)$, the quotient bundle $N=\fm/I_{r+1}$ of rank~$\mu$, and its dual that we denote by $D=N^\vee$.

Denote by $\Hilb_\eta(E)\subset\Hilb_\ud(E)$ the closure of the subvariety formed by flags of ideals with prescribed isomorphism types of the
quotient filtered algebras. Both $\Hilb_\eta(E)$ and $\Hilb_\ud(E)$ are usually singular. Consider an embedding of $\Hilb_\ud(E)$ to some compact nonsingular variety $M$ such that the action of $\GL(E)$ and the bundle~$D$ extend to~$M$.

\th[propGys] Proposition. The Thom polynomial is given by the formula
\begin{equation}\label{SigmaGys}
 \Tp_\eta=p_*(P'_\eta P''_\eta c_{\rm top}(D,F)),\qquad p:M\to\pt,
\end{equation}
where $p_*$ is the push-forward, or Gysin homomorphism associated with the projection $p$, $P'_\eta$ is the cohomology class represented by the subvariety $\Hilb_\eta(E)\subset M$, $P''_\eta$ is the normalizing factor of Definition~{\rm\ref{norm}}, and $c_{\rm top}$ is the Euler class of the corresponding bundle.

The result is actually independent of the
choice of~$M$. The traditional embedding space for the Hilbert
scheme is the Grassmann (or flag) variety of subspaces of fixed
dimension in the vector space~$\fm/\fm^{K+1}$ of jets of functions on~$E$. The main geometric idea of our approach is a construction of a new
embedding space for the Hilbert scheme that we call the
\emph{nonassociative Hilbert scheme}. The advantage of this
construction is that all ingredients of~\eqref{SigmaGys} in this
case can be identified and computed explicitly.

\medskip
{\em Proof}. Consider the variety parametrizing pairs of the form $(I,f)$ where $I\subset (\fm/\fm^{K+1})$ is an ideal of codimension~$\mu$ whose quotient algebra is isomorphic to one of the algebras defining the singularity~$\eta$, and $f$ is the $K$-jet of a map whose all components $f_i$ belong to~$I$. Denote by $\widetilde\eta$ the closure of this variety in~$\Hilb^{(\mu)}(E)\times J^K_0(E,F)$, where $\Hilb^{(\mu)}(E)$ is the \emph{punctual Hilbert scheme} parametrizing all $\mu$-codimensional ideals in $\fm/\fm^{K+1}$. The projection $p_2:\Hilb^{(\mu)}(E)\times J^K_0(E,F)\to J^K_0(E,F)$ to the second factor is proper and the details in the discussion of Definition~\ref{Ksing} show that we have the equality of classes of rational equivalence of subvarieties
\begin{equation}
\Tp_\eta\smallfrown[J^K_0(E,F)]=p_{2*}(\widetilde\eta).
 \label{p1}\end{equation}

By definition, $\widetilde\eta$ is the total space of a vector subbundle in the trivial bundle $J^K_0(E,F)\times\Hilb^{(\mu)}(E)\to \Hilb^{(\mu)}(E)$ restricted to the subvariety in $\Hilb^{(\mu)}(E)$ formed by algebras of the given isomorphism type. The Euler class of the quotient bundle of this subbundle is equal to $c_{\rm top}(\Hom(D,F))$ where $D=(\fm/I)^\vee$ is the dual of the canonical quotient bundle over the Hilbert scheme. The base of the bundle is of the form $\pi(\Hilb_\eta(E))$ where $\pi:\Hilb_{\ud}(E)\to\Hilb^{(\mu)}$ is the forgetful map that ignores the filtration on the quotient nilpotent algebra. Moreover, by the definition of the normalizing factor, we have
$$[\pi(\Hilb_{\eta}(E))]=\pi_*(P''_\eta\smallfrown[\Hilb_{\eta}(E)]).$$
Therefore,
\begin{equation}
[\widetilde\eta]=(\pi\times\Id)_*(P''_\eta c_{\rm top}(\Hom(D,F))
\smallfrown[\Hilb_{\eta}(E)\times J^K_0(E,F)]]).
 \label{p2}\end{equation}

Combining~\eqref{p1} and~\eqref{p2} we obtain finally
$$\Tp_\eta\smallfrown[J^K_0(E,F)]=
 \tilde p_{2*}(P''_\eta c_{\rm top}(\Hom(D,F))\smallfrown
   [\Hilb_\eta(E)\times J^K_0(E,F)])$$
where $\tilde p_2=p_2\circ(\pi\times\Id):\Hilb_\eta(E)\times J^K_0(E,F)\to J^K_0(E,F)$ is the projection to the second factor.

A cartesian product with a vector space induces isomorphisms both in cohomology and in the Chow groups. Therefore, the obtained equality is equivalent to the one of Proposition~\ref{propGys}.\QED

\section{The nonassociative Hilbert scheme}\label{sec6}

Denote by $\fm$ as above the maximal ideal in the ring of function
germs at the origin of a vector space~$E$ (or $K$-jets of
functions for some $K\gg0$). Let $I\subset\fm$ be an ideal of
finite codimension~$\mu$. Denote by $N=\fm/I$ the corresponding
nilpotent quotient algebra. Then we have two natural linear maps
$$\psi_1:E^\vee\to N,\qquad \psi_2:\Sym^2N\to N.$$
The first one is the restriction of the natural projection
$\fm\to\fm/I=N$ to the subspace $E^\vee\subset\fm$ consisting of
\emph{linear} functions. The second one is the structure morphism
for the multiplication law in the algebra~$N$. The following
observation is obvious: \emph{the induced linear map
$\psi_1\oplus\psi_2:E^\vee\oplus\Sym^2 N\to N$ is surjective.}

Indeed, consider a polynomial representing any element of~$N$.
This polynomial can be represented as a sum of linear terms which
are in the image of~$\psi_1$ and terms of order greater or equal
to two which are in the image of~$\psi_2$.

Conversely, let~$N$ be a $\mu$-dimensional commutative associative
nilpotent algebra with the structure morphism~$\psi_2:\Sym^2N\to
N$, and $\psi_1:E^\vee\to N$ be an arbitrary linear map such that
$\psi_1\oplus\psi_2$ is surjective. Then $\psi_1$ extends uniquely
to an algebra homomorphism $\widetilde\psi_1:\fm\to N$ for
$E^\vee$ generates~$\fm$. The homomorphism $\widetilde\psi_1$ is
surjective since $\psi_1\oplus\psi_2$ is surjective. It follows
that~$N$ can be identified with $\fm/I$ where
$I=\ker(\widetilde\psi_1)$.  We arrive at the following
conclusion.

\th Lemma. There is a one-to-one correspondence between the set of
$\mu$-codimensional ideals in~$\fm$ and the set of isomorphism
classes of pairs $(\psi_1,\psi_2)$ where $\psi_2:\Sym^2N\to N$ is
an associative commutative nilpotent algebra structure on a
$\mu$-dimensional vector space~$N$ and~$\psi_1:E^\vee\to N$ is a
linear map such that $\psi_1\oplus\psi_2$ is surjective.\QED

For the next definition we assume that the nilpotent algebra
structure that we construct on a given vector space~$N$ is
compatible with a given filtration
$$N=N_1\supset N_2\supset\dots\supset N_{r+1}=0,\qquad
            \dim N_{k}/N_{k+1}=d_k.$$
Demote by~$\widetilde M_r\subset\Hom(E^\vee\oplus \Sym^2N,N)$ the
subset determined by the following two conditions: 1) the second
component $\psi_2$ of $\psi_1\oplus\psi_2\in \widetilde M_r$
satisfies $\psi_2(N_i\*N_j)\subset N_{i+j}$ for all~$i$ and~$j$;
2) the map $\psi_1\oplus\psi_2$ is surjective. The first condition
determines a vector subspace in $\Hom(E^\vee\oplus \Sym^2N,N)$.
The second one determines a Zarisski open subset in that subspace.
The group~$\GLN$ of preserving the flag~$N_\bullet$ linear
transformations of~$N$ acts naturally on~$\widetilde M_r$.

\ex[noassdef] Definition. The \emph{nonassociative Hilbert scheme}
is the orbit space $M_r$ of the action of~$\GLN$ on $\widetilde
M_r$.

\th[freeact] Proposition. The action of~$\GLN$ on $\widetilde M_r$
is free. The variety~$M_r$ is smooth and compact.

For the proof of this proposition we provide an independent
construction of~$M_r$ which is equivalent to that used in
Definition~\ref{noassdef}. It is more convenient to use in this
construction the dual coalgebra structure on~$D=N^\vee$. Namely,
we construct~$M_r$ as the moduli space of flags
$$0=D_0\subset D_1\subset\dots\subset D_r=D,
  \qquad\dim(D_k/D_{k-1})=d_k,$$
equipped with an \emph{injective} linear map $D\to E\oplus\Sym^2D$
such that $D_k\subset E\oplus S_k$ for $k=1,\dots,r$, where
$$S_k\subset\Sym^2D_{k-1}\subset\Sym^2 D$$ is generated by subspaces of the form
$D_i\*D_j$ with $i+j\le k$.

The construction goes by induction in~$r$. In the case $r=1$ we
have $S_1=0$, therefore, $D_1$ must be embedded to~$E$. We define
$M_1=G_{d_1}(E)$, the Grassmann manifold with the tautological
rank~$d_1$ bundle~$D_1$ over it.

Assume that the variety~$M_{r-1}$ is already constructed with the
corresponding tautological flag of bundles $D_1\subset\dots\subset
D_{r-1}$ over it and with the embedding of subbundles
$D_{r-1}\subset E\oplus S_{r-1}$. Remark that~$S_r$ is determined
by $D_1,\dots,D_{r-1}$, therefore, $S_r$ can be regarded as a
bundle over $M_{r-1}$. The manifold~$M_r$ should parameterize
subspaces $D_r$ in $E\oplus S_r$ containing the subspaces
$D_{r-1}$ constructed on the previous step. According to this, we
define $M_r=G_{d_r}((E\oplus S_r)/D_{r-1})$, the corresponding
Grassmann bundle over $M_{r-1}$.

By construction, $M_r$ is defined as the total space of a tower of
fibrations with smooth compact fibers:
$$\xymatrix{
 M_r\ar[rr]^{G_{d_r}(E_r/D_{r-1})}&&
 M_{r-1}\ar[rr]^{G_{d_{r-1}}(E_{r-1}/D_{r-2})}&&
 \dots\ar[rr]^{G_{1}(E)}&&\pt
},\qquad E_k=E\oplus S_k.$$
This implies
Proposition~\ref{freeact}.\QED

The just constructed manifold~$M_r$ is equipped with the diagram
of bundles and subbundles
$$\xymatrix{
  D_1\ar@{_{(}->}[d]\ar@{^{(}->}[r]&
  D_2\ar@{_{(}->}[d]\ar@{^{(}->}[r]&\cdots\ar@{^{(}->}[r]&
  D_r\ar@{_{(}->}[d]\\
  E=E_1\ar@{^{(}->}[r]&E_2\ar@{^{(}->}[r]&
     \cdots\ar@{^{(}->}[r]&E_r,&&E_k=E\oplus S_k.}$$
The projection of the subbundle $D_r\subset E\oplus S_k$ to the
summands $E$ and $S_r$ determines a linear map $D_r\to E$ and a
canonical filtered commutative coalgebra structure $D_r\to S_r$ on
the fibers of~$D_r$. The dual bundle $N_r=D_r^\vee$ equipped,
respectively, with a surjective bundle map $E^\vee\oplus
S_r^\vee\to N_r$ defining a canonical linear map $E^\vee\to N_r$
and a canonical commutative filtered algebra structure on the
fibers of~$N_r$. The following assertion is a reformulation of
Lemma~\ref{noassdef}.

\th[Hilb] Proposition. The Hilbert scheme~$\Hilb_\ud(E)$ is
isomorphic to the sublocus in $M_r$ consisting of points for which
the canonical commutative filtered algebra on the corresponding
fiber of~$N_r$ is associative. Its subvariety~$\Hilb_\eta(E)$
consists of the points of~$M_r$ for which the canonical algebra on
the corresponding fiber of~$N_r$ is isomorphic to one of the
algebras defining the singularity type~$\eta$.\QED

\section{The final computation of the Thom polynomial}\label{sec7}

We are able now to complete the computation of Thom polynomial
applying Proposition~\ref{propGys} to the projection
\begin{equation}\label{pmap}
 p:M_r\to\pt
\end{equation}
where $M_r$ is the nonassociative Hilbert scheme constructed in
the previous section. By proposition~\ref{Hilb}, the polynomial $P_\eta'(t_1,\dots,t_\mu)$ representing Poincar\'e
dual of the locus $\Hilb_\eta(E)\subset M_r$ is the local
characteristic invariant of
Definition~\ref{lochi} expressed in terms of the Chern roots
$-t_1,\dots,-t_\mu$ of the bundles $D_1/D_0$, \dots,
$D_{r}/D_{r-1}$ over $M_r$. It remains to compute the
homomorphism~$p_*$. Various alternative methods are known for the
computation of the Gysin homomorphism. We prefer to use the
machinery developed in~\cite{Kaz11} by two reasons. First, with this
approach the formulas become quite simple. Second, this approach
allows one to keep the Chern classes $c(E)$ and $c(F)$ in the
relative combination $c(F{-}E)$ in all intermediate steps of
computations. Here is the basic relation from~\cite{Kaz11}.

\th[pstar] Proposition. Let $E$ and $F$ be vector bundles of
relative rank $\ell=\rk F-\rk E$ defined over a smooth base $X$.
Denote by $D$ the canonical rank $d$ bundle over the total space of
the Grassmann fibration $q:G_d(E)\to X$. Then the Gysin
homomorphism $q_*$ acts by the following formula:
$$q_*:P(t_1,\dots,t_d)\;c_{\rm top}(\Hom(D,F))\mapsto
  P(t_1,\dots,t_d)\;t_1^{\ell+d}t_2^{\ell+d-1}\dots t_d^{\ell+1}
   \prod_{1\le i<j\le d}\!\!\!(t_i-t_j).$$

The meaning of this relation is the following. The class $c_{\rm
top}$ is the Euler class of the corresponding bundle. The
polynomial~$P$ on the left hand side represents any polynomial
combination of the Chern classes $c_i(D)$ expressed formally as a
symmetric function in the corresponding Chern roots
$-t_1,\dots,-t_d$, i.e. we set formally $c(D)=\prod_{i=1}^d(1-t_i)$. \emph{The meaning of the same variables $t_k$
on the right hand side is different. Namely, the Laurent expansion of the right hand side is the generating series for the coefficients of the polylinear expansion in the Chern classes $c_i(F{-}E)$.} Informally, the monomial $t_k^i$
in the monomial expansion of the right hand side represents the class $c_i(F{-}E)$ (for any~$k$).

In order to apply this proposition we represent the mapping~$p$
in~\eqref{pmap} as the composition
$$M_r\too[q_r]M_{r-1}\too[q_{r-1}]\dots\too[q_2] M_1\too[q_1] M$$
and decompose also
$$c_{\rm top}(\Hom(D_r,F))
 =c_{\rm top}(\Hom(D_{r-1},F))\,c_{\rm top}(\Hom(D_r/D_{r-1},F))
 =\prod_{k=1}^rc_{\rm top}(\Hom(D_k/D_{k-1},F)).$$
Recall that $q_r$ is the Grassmann fibration
$G_{d_r}(E_r/D_{r-1})$ over~$M_{r-1}$ where $E_r=E\oplus S_r$. By
Proposition~\ref{pstar}, in order to apply $q_{r*}$ we need to
replace the factor $c_{\rm top}(\Hom(D_r/D_{r-1}))$ by the factor
\begin{equation}\label{eqst1}
s_1^{\ell_r+d_r}s_2^{\ell_r+d_r-1}\dots s_{d_r}^{\ell_r+1}
   \prod_{1\le i<j\le d_r}\!\!\!(s_i-s_j)
\end{equation}
where we denote temporarily $s_i=t_{d_1+\dots+d_{r-1}+i}$ and
$\ell_r=\rk F-\rk E_r/D_{r-1}=\ell+\rk D_{r-1}-\rk S_r$. The
product of~\eqref{eqst1} with $P_\eta\,c_{\rm top}(\Hom(D_{r-1},F))$
represents the class $q_{r*}(P_\eta\,c_{\rm top}(\Hom(D_r,F)))$
where~$s_m^k$ in the monomial expansion denotes the class
$$c_k(F{-}E_r/D_{r-1})=c_k(F{-}E+D_{r-1}{-}S_{r}).$$
We prefer to have an expression where $s_m^k$ would represent the
class $c_k=c_k(F{-}E)$. By the Whitney formula we need to replace
every monomial $s_m^k$ by the expression
$$s_m^k\leadsto s_m^k\sum_{j\ge0}
 c_j(D_{r-1}-S_r)s_m^{-j}
 =s_m^k\frac{\prod_{i=1}^{d_1+\dots+d_{r-1}}(1-t_i/s_m)}{\prod_\l(1-t_\l/s_m)}$$
where $\l$ runs over the list of Chern roots of the bundle~$S_r$
expressed in terms of the Chern roots of $D_k$: $t_\l=t_i+t_j$,
$1\le i\le j$, $w(i)+w(j)\le r$. We see that passing from the
presentation in terms of the Chern classes $c_k(F-E_r/D_{r-1})$ to
the Chern classes in $c_k=c_k(F{-}E)$ is accounted by a factor.
Multiplying all these factors for $m=1,\dots, d_r$ and multiplying
by~\eqref{eqst1} we obtain that the actual action of $q_{r*}$ is
correctly represented by the factor
$$\frac{\displaystyle
  \prod_{i=1}^{d_r}s_i^{\ell+1+d_r-i}\;\;
  \prod_{j=1}^{d_r}\prod_{i=1}^{d_1+\dots+d_{r-1}\!\!\!\!\!\!}\!\!\!\!\!(s_j-t_i)
  \prod_{1\le i<j\le d_r}\!\!\!(s_j-s_i)}
  {\displaystyle\prod_{k=1}^{d_r}\prod_\l(s_k-t_\l)}.$$
where $s_i$ still denotes $t_{d_1+\dots+d_{r-1}+i}$. The next
homomorphism $(q_{r-1})_*$ introduces another similar factor, and
so on. Multiplying all these factors corresponding to $q_{r*}$,
\dots, $q_{1*}$ we arrive at the function of Theorem~\ref{th1}, up
to the factor $\prod_{i=1}^\mu t_i^{\ell+1}$ which is accounted by
the shift of the indices of the Chern classes $c_i$ in
Eq.~\ref{Thser}. Theorem~\ref{th1} is proved.\QED

\section{Examples and computations}\label{sec8}

In this Section we present some explicit computations based on the formula of Theorem~\ref{th1}. We use traditional notations for singularity types corresponding to isomorphism classes of nilpotent algebras, see~\cite{DS6}. An algebra can be characterized by a collection of generators of the ideal whose quotient is isomorphic to the given algebra. However, it is more demonstrative to use a different convention identifying algebras by their monomial bases. Namely, let a list of monomials in triangular brackets denote the quotient algebra over the \emph{monomial} ideal spanned by the monomials which are \emph{not} present in the list. Besides, we denote by $\fm_{x,y,\dots}$ the maximal ideal in the local ring $\C[[x,y\dots]]$ of infinite power series in the variables $x,y,\dots$. The discussion of this section involves the following nilpotent algebras:
{\renewcommand{\arraystretch}{1.5}
\begin{longtable}[c]{|c|c|c|}\hline
Notation&Algebra&dimension\\\hline
    $\S^{\mu}$&$\displaystyle
  \<x_1,\dots,x_\mu\>=\frac{\fm_{x_1,\dots,x_\mu}}{\fm^2}$&$\mu$\\\hline
      $A_\mu$&$\displaystyle
  \<x,x^2,\dots,x^\mu\>=\frac{\fm_x}{(x^{\mu+1})}$&$\mu$\\\hline
    $I_{a,b}$&$\displaystyle
  \frac{\<x,x^2,\dots,x^a,y,y^2,\dots,y^b\>}
   {(x^a+x^b)}=\frac{\fm_{x,y}}{(xy,x^a+y^b)}$&$a+b-1$\\\hline
  $III_{a,b}$&$\displaystyle
  \<x,x^2,\dots,x^{a-1},y,y^2,\dots,y^{b-1}\>
  =\frac{\fm_{x,y}}{(xy,x^a,    y^b)}$&$a+b-2$\\\hline
 $\Phi_{m,r}$&$\displaystyle
  \frac{\<x_i,y_j,y_k^2\>}
  {(y_j^2-y_{k}^2)},\quad
  1\le i\le r,~1\le j,k\le m-r$&$m+1$\\\hline
   $\S^{a,b}$&$\displaystyle
  \<x_i,y_k,x_ix_j,x_iy_k\>,
   \quad 1\le i,j\le b,~1\le k\le a-b$&$\mu(a,b)$\\\hline
   $\S^{j_1,\dots,j_m}$&$\displaystyle
  \<x_{i_r}\!\!\cdots x_{i_1}\>,
   \quad 1\le r\le m,~i_1\ge\dots\ge i_r\ge1,$&$\mu(j_1,\dots,j_m)$\\
   $j_1\ge\dots\ge j_m\ge1$&
    $1\le i_1\le j_1,~\dots,~1\le i_r\le j_r$&\\\hline
   $C_\l$&$\displaystyle
  \frac{\<x,y,z,x^2,y^2,z^2,xy,yz,zx\>}
  {(y^2+2xz,2yz,x^2+6xz+\g z^2)}$&$6$\\\hline
\end{longtable}}

The dimension $\mu(j_1,\dots,j_m)$ of the Thom-Boardman algebra $\S^{j_1,\dots,j_m}$ is equal to the number of monomials in the basis from its description above. In particular, for $m=2$ we have
$$\mu(a,b)=a+\frac{b(b+1)}{2}+b(a-b)=a(b+1)-\frac{b(b-1)}{2}.$$

\subsection{Thom-Porteous formula}

The locus of the Thom-Porteous singularity $\S^\mu$ for a given map $f:X\to Y$ is the locus of points $x\in X$ characterized by the condition that the derivative $df:T_xX\to T_{f(x)}Y$ has kernel rank at least~$\mu$. The nilpotent algebra of this singularity has trivial multiplication. Therefore, it admits a natural filtration with the dimension vector~$\ud=(\mu)$ consisting of only one entry. By Corollary~\ref{cor1}, the generating series for the Schur basis expansion of the Thom polynomial is just the monomial
$$\tK_\mu(t_1,\dots,t_\mu)=t_1^{\mu-1}t_2^{\mu-1}\dots t_\mu^{\mu-1}$$
which recovers the Porteous formula
$$\Tp_{\S^\mu}=\Delta_{\ell+\mu,\dots,\ell+\mu}.$$

Remark that not only the answer for the singularity $\S^\mu$ agrees with the Porteous formula but also the resolution used in the proof for this particular case coincides with the `standard' one.

The Thom-Porteous singularity is the only one for which the generating series is a polynomial. In all other cases the rational  generating series has a non-trivial denominator leading to an infinite monomial expansion.

\subsection{Example: Thom series for $I_{2,2}$ and $III_{2,3}$}

Consider the $3$-dimensional vector space~$N$ with the basis $e_1,e_2,e_3$ and an algebra structure on it given by the multiplication law
$$e_1\cdot e_1=q_{11}^3e_3,\qquad
  e_1\cdot e_2=e_2\cdot e_1=q_{12}^3e_3,\qquad
  e_2\cdot e_2=q_{2,2}^3e_3$$
and all other products of basic vectors being trivial. The complex coefficients $q_{11}^3,q_{12}^3,q_{22}^3$ form the $3$-dimensional vector space $\Alg(2,1)$. The multiplication defines a quadratic form on $N^{(1)}=N/N_2\simeq\C^2$ taking values in $N_2\simeq\C$. If this form is nondegenerate, then the algebra is isomorphic to~$I_{2,2}=\frac{\fm_{x,y}}{(x^2,y^2)}$. It follows that the Thom polynomial for the singularity $I_{2,2}$ admits the generating series
\begin{align*}
I_{2,2}:\qquad K_{2,1}&=\frac{t_1(t_2-t_1)(t_3-t_1)(t_3-t_2)}
  {(t_3-2t_1)(t_3-t_1-t_2)(t_3-2t_2)}\\
  &=t_1^1t_2^1t_3^{-1}-t_1^2t_2^0t_3^{-1}-2t_1^3t_2^0t_3^{-2}
   +2t_1^1t_2^2t_3^{-2}-4t_1^4t_2^0t_3^{-3}+t_1^2t_2^2t_3^{-3}\\
   &\quad+(4t_1^1t_2^3t_3^{-3}-t_1^3t_2^1t_3^{-3})-8t_1^5t_2^0t_3^{-4}+3t_1^2t_2^3t_3^{-4}
   +(8t_1^1t_2^4t_3^{-4}-3t_1^4t_2^1t_3^{-4})+\dots
\end{align*}
leading to the Thom polynomial
\begin{align*}
\Tp_{I_{2,2}}&=c_{\ell+2}^2c_\ell - c_{\ell+3}c_{\ell+1}c_\ell
 -2c_{\ell+4}c_{\ell+1}c_{\ell-1}+2c_{\ell+3}c_{\ell+2}c_{\ell-1}
 -4c_{\ell+5}c_{\ell+1}c_{\ell-2}+c_{\ell+3}^2c_{\ell-2}\\
 &\qquad+3c_{\ell+4}c_{\ell+2}c_{\ell-2}-8c_{\ell+6}c_{\ell+1}c_{\ell-3}
 +3c_{\ell+4}c_{\ell+3}c_{\ell-3}+5c_{\ell+5}c_{\ell+2}c_{\ell-3}+\dots
 \end{align*}
where the dots denote the terms vanishing for $\ell>3$. We see that different monomials of the generating series may contribute to one and the same monomial in the Chern classes. A consequence of this phenomenon is that the generating series is not unique. In particular, another choice of the filtration can lead to a different generating series for the same Thom polynomial. Consider, for example, the filtration with the dimension vector $(0,1,1,0,1)$. The generic algebra with this filtration is also isomorphic to~$I_{2,2}$. One should take into account that such a filtration on $I_{2,2}$ is \emph{not} natural in the sense of Definition~\ref{fnat}: the choice of the subspace $N_3\subset N$ corresponds to the choice of one of the two zero lines of the quadratic form mentioned above. By Theorem~\ref{th1}, we obtain an alternative generating series for the Thom polynomial~$\Tp_{I_{2,2}}$:
$$I_{2,2}:\qquad \frac12 K_{0,1,1,0,1}=\frac{(t_2-t_1)(t_3-t_1)(t_3-t_2)}
  {2(t_3-2t_1)(t_3-t_1-t_2)}.$$

The fact that the two generating functions determine the same Thom polynomial can be checked algebraically. Namely, represent the difference of these two functions in the form
$$\frac12K_{0,1,1,0,1}-K_{2,1}=\frac{(t_3-2t_1-2t_2)(t_2-t_1)(t_3-t_2)(t_3-t_1)}
 {2(t_3-2t_1)(t_3-t_1-t_2)(t_3-2t_2)}.$$
This function is skew-symmetric with respect to $t_1$ and $t_2$. More precisely, the coefficients of its expansion as a Laurent series in $t_3^{-1}$ are skew-symmetric Lauret polynomials in $t_1$ and $t_2$. Therefore, its contribution to each monomial in $c$-variables cancel.

Corollary~\ref{cor1} provides also a representation of the Thom polynomial in the Schur basis: expanding the function
\begin{align*}
\frac12 \tK_{0,1,1,0,1}&=\frac{t_2t_3^2}
  {2(t_3-2t_1)(t_3-t_1-t_2)}
  =\sum_{k,j=0}^\infty
  \Bigl(\sum_{i=0}^k\binom{i+j}{i}2^{k-i-1}\Bigr)
  t_1^kt_2^{j+1}t_3^{-k-j}
\end{align*}
we get
\begin{multline*}
\Tp_{I_{2,2}}
  =\sum_{k,j=0}^\infty
  \Bigl(\sum_{i=0}^k\binom{i+j}{i}2^{k-i-1}\Bigr)
  \Delta_{\ell+1+k,\ell+2+j,\ell+1-i-j}
 =\Delta_{\ell+2,\ell+2,\ell}+3\Delta_{\ell+3,\ell+2,\ell-1}\\
 +3\Delta_{\ell+3,\ell+3,\ell-2}+7\Delta_{\ell+4,\ell+2,\ell-2}
 +10\Delta_{\ell+4,\ell+3,\ell-3}+15\Delta_{\ell+5,\ell+2,\ell-3}+\dots.
 \end{multline*}
It is easy to check applying determinantal formulas for Schur functions that the two obtained expansions for the Thom polynomial of the singularity~$I_{2,2}$ agree.

\medskip
Let us turn back to the space $\Alg(2,1)$. The discriminant has the equation $q_{11}^3q_{2,2}^3-(q_{1,2}^3)^2=0$ of multidegree $2(t_3-t_2-t_1)$. The algebras from this discriminant are isomorphic to $III_{2,3}$. This leads to the following generating series for its Thom polynomial:
$$III_{2,3}:\qquad 2(t_3-t_2-t_1)K_{2,1}
 =\frac{2t_1(t_3-t_2-t_1)(t_2-t_1)(t_3-t_1)(t_3-t_2)}
  {(t_3-2t_1)(t_3-t_1-t_2)(t_3-2t_2)}.$$

It is interesting to observe that this algebra admits two more natural filtrations, those with the dimension vectors $(1,2)$ and $(0,1,1,1)$. Moreover, the algebras isomorphic to $III_{2,3}$ form Zariski open subsets in both~$\Alg(1,2)$ and $\Alg(0,1,1,1)$. This leads to two other possible choices of the generating series for the singularity~$III_{2,3}$:
$$III_{2,3}:\qquad  K_{1,2}=\frac{t_2(t_2-t_1)(t_3-t_1)(t_3-t_2)}
  {(t_3-2t_1)(t_2-2t_1)},\qquad
 K_{0,1,1,1}=\frac{(t_2-t_1)(t_3-t_1)(t_3-t_2)}{t_3-2t_1}.$$
The last function is certainly the simplest one among the three possible choices of the generating series for the $III_{2,3}$ singularity presented above. It provides a particulary simple expression for the Thom polynomial in the Schur basis: expanding the series
$$\tK_{0,1,1,1}=\frac{t_2t_3^2}{t_3-2t_1}=t_2t_3+2t_1t_2+4t_1^2t_2t_3^{-1}
 +8t_1^3t_2t_3^{-2}+\dots$$
 we get (cf.~\cite{FRloc})
$$\Tp_{III_{2,3}}=\sum_{i=1}^{\ell+2}
      2^i\Delta_{\ell+1+i,\ell+2,\ell+2-i}.$$

\subsection{Summary on the computed Thom polynomials}\label{sectab}

An analysis of algebras corresponding to different choices of filtrations with relatively small number of nonzero entries in the dimension vector, similar to the one of the previous section, leads to the computation of generating series for Thom polynomials of particular classes of singularities. A part of these computations is presented in the tables below. In these tables, we denote by~$d$ the homogeneous degree of the generating function. It is related to the codimension of the singularity by the equality
$$\codim\eta=(\ell+1)\,\mu+d.$$
One may treat~$d$ as the `virtual codimension' of the singularity~$\eta$ corresponding to the case~$\ell=-1$. The number~$c$ in the last column is the degree of the denominator (the number of its linear factors) that measures the `complexity' of the corresponding rational function from the point of view of its power series expansion.

\def\dow#1{\raisebox{-1.5ex}[0pt][0pt]{#1}}
\renewcommand{\arraystretch}{1.3}
\begin{longtable}[c]{|c|c|c|l|c|c|c|}
\hline $\eta$&$\mu$&$d$&\qquad\qquad Generating series for $\Tp_\eta$&$c$\\\hline\hline
\endhead
\hline\endfoot
 $A_1=\S^1$&$1$&$0$&$K_{1}$&0\\
\hline
 $A_2=\S^{1,1}$&2&0&$K_{1,1}$&1\\
\hline
 $\S^2$&2&2&$K_{2}$&0\\
\hline
 $A_3=\S^{1,1,1}$&3&0&$K_{1,1,1}$&3\\
\hline
 \dow{$\Phi_{2,0}=I_{2,2}$}&\dow3&\dow1&$K_{2,1}$&3\\*
    &&&$K_{0,1,1,0,1}\cdot\frac12$&2\\
\hline
 &&&$K_{2,1}\cdot 2(t_3-t_1-t_2)$&2\\*
 $\Phi_{2,1}{=}III_{2,3}$&3&2&$K_{1,2}$&2\\*
 &&&$K_{0,1,1,1}$&1\\
\hline
 $\S^3$&3&6&$K_{3}$&0\\
\hline\hline
 $A_4{=}\S^{1,1,1,1}$&4&0&$K_{1,1,1,1}\cdot(t_4-2t_1-t_2)$&7\\
\hline
 \dow{$I_{2,3}$}&\dow4&\dow1&$K_{2,1,1}\cdot
  ((t_4{-}2t_1{-}t_2)(t_4{-}t_1{-}2t_2)
     -(t_4{-}t_2{-}t_3)(t_4{-}t_1{-}t_3))$&8\\*
 &&&$K_{0,1,1,1,0,1}$&5\\
\hline
 \dow{$III_{3,3}$}&\dow4&\dow2&$K_{2,2}$&6\\*
 &&&$K_{0,0,1,1,0,1,0,1}\cdot\frac12$&4\\
\hline
 &&&$K_{1,2,1}$&5\\*
 \dow{$III_{2,4}$}&\dow4&\dow2&$K_{1,1,2}$&5\\*
 &&&$K_{0,1,0,1,1,1}$&4\\*
 &&&$K_{0,0,1,0,1,1,0,0,1}$&4\\
\hline
 &&&$K_{2,2}\cdot 2(t_3+t_4-2t_1-2t_2)$&6\\*
 $\S^{2,1}$&4&3&$K_{0,1,1,0,2}$&4\\*
 &&&$K_{0,1,1,1,1}$&3\\
\hline
 \dow{$\Phi_{3,0}$}&\dow4&\dow3&$K_{3,1}$&6\\*
 &&&$K_{0,0,1,1,1,0,0,1}\cdot \frac{-1}{4}(t_4-2t_1)$&3\\
\hline
 &&&$K_{3,1}\cdot(3t_4-2t_1-2t_2-2t_3)$&6\\*
 $\Phi_{3,1}$&4&4&$K_{0,2,1,1}$&3\\*
 &&&$K_{0,0,1,1,1,0,1}\cdot \frac{1}{2}$&2\\
\hline
 \dow{$\Phi_{3,2}$}&\dow4&\dow6&$K_{3,1}\cdot4(t_4{-}t_1{-}t_2)(t_4{-}t_2{-}t_3)(t_4{-}t_1{-}t_3)$&3\\*
 &&&$K_{0,1,2,1}$&1\\
\hline
 $\S^4$&4&12&$K_{4}$&0\\
 \hline\hline
 \dow{$\S^{2,1,1}$}&\dow6&\dow4&$K_{0,1,1,1,1,1,1}\cdot(t_6-2t_1-t_2)(t_5-2t_2)$&12\\*
 &&&\quad${}=K_{0,0,1,0,1,1,0,1,1,0,1}\cdot(t_6-2t_1-t_2)$&\\
\hline
 $\S^{2,2,1}$&8&7&$K_{0, 0, 1, 1, 0, 1, 1, 1, 1, 1, 1} \cdot
 (t_7{-}2t_1{-}t_2)(t_8{-}2t_1{-}t_2)(t_8{-}t_1{-}2t_2)$&24\\
\hline
 \dow{$\S^{2,2,2}$}&\dow9&\dow9&$K_{2,3,4}\cdot\prod_{k=6}^9((t_k-2t_1-t_2)(t_k-t_1-2t_2))$&45\\*
  &&&$K_{0, 0,0, 1, 1, 0,0, 1, 1, 1,0,1, 1, 1, 1} \cdot
  t_1\prod_{k=7}^9(t_k{-}2t_1{-}t_2)\prod_{k=8}^9(t_k{-}t_1{-}2t_2)$&33\\
\hline
  &&&$K_{3,3}\cdot4(t_4+t_5+t_6-2t_1-2t_2-2t_3)$&18\\*
  &&&$K_{0,2,1,0,3}$&15\\*
  $C_\g$&6&4&$K_{0,0,0,1,1,1,0,0,1,1,1}\cdot\frac12$&11\\*
  &&&$K_{0,0,1,2,0,1,0,2}$&13\\*
  &&&$K_{0,0,0,0,1,1,1,0,0,1,0,1,0,1}\cdot\frac12$&11\\
\end{longtable}

\renewcommand{\arraystretch}{2}
\begin{longtable}[c]{|c|c|c|l|c|c|c|}
\hline $\eta$&$\mu$&$d$&Generating series&$c$\\\hline\hline
 $\S^r$&$r$&$r(r-1)$&$K_{r}$&0\\
\hline
 $\Phi_{m,r}$&$m+1$&$\binom{m}{2}{+}\binom{r+1}{2}$&$K_{0,m-r,r,1}$&$\binom{m-r+1}{2}$\\
\hline
 $\S^{a,b}$&$a(b{+}1){-}\binom{b}{2}$&$(a{-}1)\mu{-}(a{-}b)b$
 &$K_{0,b,a-b,\binom{b+1}{2},b(a-b)}$
 &$\binom{b+1}{2}^2{+}b(a{-}b)\bigl(b(a{-}b)+\binom{b+1}2\bigr)$\\
\hline
\end{longtable}

Every line of the presented tables can be regarded as an independent theorem. We will not provide complete proofs leaving most of the (quite elementary) checks to the reader. In the following sections we discuss the details of some most interesting cases manifesting particular phenomena appearing in the course of computations. But here we finish with just one useful observation.

Define the grading in the ring $\C[[x_1,\dots,x_m]]$ by assigning certain integer weights $\deg x_i=\o_i$ to the variables. For an integer $K>0$, denote by $I_{\o_1,\dots,\o_m}^{(K)}$ the ideal generated by all monomials of degree greater than~$K$ and set $N_{\o_1,\dots,\o_m}^{(K)} =\frac{\fm_{x_1,\dots,x_m}}{I_{\o_1,\dots,\o_m}^{(K)}}$. This nilpotent algebra is graded, and we denote by  $\ud_{\o_1,\dots,\o_m}^{(K)}$ the dimension vector of the associated filtration.

\th[propquasi] Proposition. Algebras isomorphic to $N_{\o_1,\dots,\o_m}^{(K)}$ form a Zariski open subset in the subvariety of $\Alg(\ud_{\o_1,\dots,\o_m}^{(K)})$ consisting of associative algebras. If the associativity equation is empty, then a generic algebra in $\Alg(\ud_{\o_1,\dots,\o_m}^{(K)})$ is isomorphic to $N_{\o_1,\dots,\o_m}^{(K)}$.

Indeed, consider a structure of a commutative associative filtered algebra on a given vector space~$N$. Chose generic vectors in the subspaces of filtrations $\o_1,\dots,\o_m$ and denote these vectors by $x_1,\dots,x_m$. Then monomials in these elements of quasihomogeneous degree smaller or equal to~$K$ are generically linearly independent and span the whole~$N$ while any monomial of degree greater than~$K$ is equal to zero. It shows that the obtained algebra is isomorphic to~$N_{\o_1,\dots,\o_m}^{(K)}$.\QED

For a given isomorphism type of algebras we try to find a quasihomogeneous filtration to represent this algebra in the form~$N_{\o_1,\dots,\o_m}^{(K)}$ or as close to this form as possible. Thus, the funny dimension vectors with many zeroes and ones in the tables of this section correspond to different choices of quasihomogeneity weights. For example,
%
\begin{align*}
 \S^\mu&\simeq N_{1_\mu}^{(1)},&\ud&=(\mu),\\
 A_\mu&\simeq N_1^{(\mu)},&\ud&=(1_\mu),\\
 III_{2,k+1}&\simeq N_{1,k}^{(k)},&\ud&=(1_{k-1},2),\\
 III_{2,k+1}&\simeq N_{2,2k-1}^{(2k)},&\ud&=(0,1,0,1,\dots,0,1,1,1),\\
 \Phi_{m,m-1}&\simeq N_{2,3_{m-1}}^{(4)},&\ud&=(0,1,m-1,1),\\
 \S^{a,b}&\simeq N_{2_b,3_{a-b}}^{(5)},
       &\ud&=(0,b,a-b,\textstyle\binom{b+1}{2},b(a-b)),\\
 \S^{2,2,2}&\simeq N_{1,1}^{((2)},&\ud&=(2,3,4),\\
 \S^{2,2,2}&\simeq N_{4,5}^{(15)},
       &\ud&=(0,0,0,1,1,0,0,1,1,1,0,1,1,1,1),\\
 \S^{2,2,1}&\simeq N_{3,4}^{(11)},
       &\ud&=(0,0,1,1,0,1,1,1,1,1,1).
\end{align*}

\subsection{Algebras of dimension at most~$4$}

The table of the previous section contains the computed generating functions for Thom polynomials for all nilpotent algebras of dimension~$\le4$. In fact, these Thom polynomials are computed in~\cite{FRloc} by a different method (partially announced). In most cases we present several alternative forms of the generating series corresponding to different choices of possible filtrations. The classification of nilpotent algebras of small dimension is finite: there is a unique algebra $\S^1=A_1$ of dimension~$1$, two algebras,~$A_2=\S^{1,1}$ and $\S^2$ of dimension~$2$, and also~$4$ isomorphism classes of algebras of dimension~$3$, and~$9$ classes of $4$-dimensional algebras whose adjacencies are presented in the following diagram.
$$\xymatrix@R=1pc{ d:&0&1&2&3&4&6&12\\
 &A_3&I_{2,2}\ar[l]&III_{2,3}\ar[l]&&&\S^{3}\ar[lll]\\
 &A_4&I_{2,3}\ar[l]&III_{3,3}\ar[l]&\Phi_{3,0}\ar@/_1pc/[ll]&\Phi_{3,1}\ar[l]\ar[ld]&
 \Phi_{3,2}\ar[l]&\S^4\ar[l]\\
 &&&III_{2,4}\ar[lu]&\S^{2,1}\ar[lu]\ar[l]}$$

These ajacencies are used in the classification of orbits in the spaces~$\Alg(\ud)$ with $\mu=\sum d_i\le 4$. It proves out, however, that the detailed study of orbits is redundant: in all cases except for $A_4$ an algebra of dimension~$\le4$ admits a natural filtration such that the algebra is isomorphic to a generic algebra with this filtration. Therefore, the generating function for its Thom polynomial can be chosen in the form~$K_{\ud}$ for a suitable dimension vector~$\ud$.

In the case of the algebra~$A_4$ there is a unique choice of the filtration. The dimension vector of this filtration is $(1,1,1,1)$. A generic algebra with this filtration is isomorphic to~$A_4$ \emph{provided that it is associative}. This is the simplest case of a non-empty associativity condition.
The associativity equation determines a hypersurface of multidegree $t_4-2t_1-t_2$, see Example~\ref{exA4}. Therefore, we obtain a generating function for $\Tp_{A_4}$ in the form $K_{1,1,1,1}(t_4-2t_1-t_2)$.

\subsection{The singularities $A_\mu$}

The Morin singularity is the one with the nilpotent algebra $\frac{\fm_x}{(x^{\mu+1})}$.  It is denoted by $A_\mu=\S^{1_\mu}=\S^{1,\dots,1}$ ($\mu$ units). A natural filtration on this algebra is induced by the grading with $\deg x=1$. The dimension vector of this filtration is $(1_\mu)=(1,\dots,1)$. It follows that the generating series for the Thom polynomial can be chosen in the form
$$A_\mu:\qquad K_{1_\mu}\cdot R_\mu$$
where $R_\mu=P_{A_\mu}$ is the multidegree of the closure of algebras isomorphic to~$A_\mu$ in $\Alg(1_\mu)$. This presentation of the Thom polynomial for the singularity~$A_\mu$ is essentially equivalent to the one found in~\cite{BS} by a different approach. The degree of~$R_\mu$ (that is, the codimension of the corresponding subvariety in~$\Alg(1_\mu)$) is shown in the following table.
$$
\renewcommand{\arraystretch}{1.3}
\begin{array}{|c|c|c|c|c|c|c|}
\hline
\mu&1&2&3&4&5&6\\\hline
\dim\Alg(1_\mu)&0&1&3&7&13&22\\\hline
\deg R_\mu&0&0&0&1&3&7\\\hline
\end{array}
$$

For $\mu\le3$ the algebras isomorphic to $A_\mu$ form a Zariski open subset in~$\Alg(1_\mu)$. Therefore, one has $R_\mu=1$ for $\mu=1,2,3$. In the case $\mu=4$ the locus of $A_4$ in $\Alg(1,1,1,1)$ is the hypersurface determined by the associativity equation, as it is explained in Example~\ref{exA4}.

For $\mu=5$ the the locus~$A_5$ is also determined by the associativity equations. This locus has codimension~$3$ but it is not a complete intersection, its ideal is generated by~$4$ (dependent) equations. They can be represented as the rank condition on a matrix whose entries are  certain coordinates in~$\Alg(1_5)$. Thus, applying the Porteous formulas one computes the multidegree~$R_5$ of this variety, see~\cite{BS}.

In the case $\mu=6$ a new phenomenon occurs. The subvariety determined by the associativity equations in~$\Alg(1_6)$ is proved to be reducible. One of its components is the closure of the set of algebras isomorphic to~$A_6$. But there is another component, of the same dimension, that corresponds to the `nets of conics' singularity. The multidegree of the component~$A_6$ is computed by means of computer algebra in~\cite{BS}.

For greater values of~$\mu$ the situation is even more complicated. The variety of solutions of the associativity equations has many components, of even different dimensions, and there is no known regular method to single out the component~$A_\mu$ and to compute its multidegree in a closed form.

\subsection{The singularities $\Phi_{m,r}$}

The singularities $\Phi_{m,r}$, $0\le r\le m-1$, are introduced and their Thom polynomials are computed in~\cite{FRloc}. In our notations, they classify orbits in $\Alg(m,1)$. Consider a filtered nilpotent algebra $N$ of dimension $m+1$ with the dimension vector $(m,1)$. The multiplication defines a quadratic form on $N^{(1)}=N/N_2\simeq\C^m$ with values in~$N_2\simeq \C$. We say that the algebra is of type $\Phi_{m,r}$ if this quadratic form is of rank~$m-r$, that is, it has a kernel of dimension~$r$.

It follows from the definition that the Thom series for this singularity can be chosen in the form
$$\Phi_{m,r}:\qquad K_{m,1} Z_r$$
where $Z_r$ is the multidegree of the variety of quadratic forms of kernel rank~$r$. An explicit form of the polynomial~$Z_r$ is computed in the framework of the theory of symmetric degeneracy loci, see~\cite{sym,HT,Pr91,Ful1,FP98}. On the other hand, there is a simpler and even more efficient way to compute the Thom polynomial in this case. Namely, the existence of an $r$-dimensional distinguished subspace (the kernel of the quadratic form) allows one to use a more detailed filtration in the algebra, the one with the dimension vector $(0,m-r,r,1)$. A generic algebra with this dimension vector is isomorphic to~$\Phi_{m,r}$. That means that $K_{0,m-r,r,1}$ is a generating function for the Thom polynomial of this singularity.

\subsection{The singularity~$\S^{a,b}$}

The Thom-Boardman singularity $\S^{a,b}$, $1\le b\le a$, is, by definition, the one associated with the quotient algebra of the coordinate ring $\C[[x_1,\dots,x_{b},y_1,\dots,y_{a-b}]]$ over the ideal spanned by the monomials of the form $y_iy_j$ and also all monomials of degree greater than~$2$. A general algorithm for the computation of the Thom polynomial for these singularities was derived in~\cite{Ro}. A closed formula for the Thom series was found in~\cite{FK06} for particular cases and in~\cite{Kaz11} for general one. In the present approach the answer is obtained immediately with essentially no computations involved.

Consider a grading on the coordinate ring assigning the degrees of the variables to be $\deg x_i=2$, $\deg y_j=3$. Then the ideal is spanned by monomials of degree~$\ge 6$. The filtration on the algebra $\S^{a,b}$ induced by this grading has the dimension vector~$(0,b,a-b,\binom{b+1}{2},b(a-b))$. Moreover, a generic algebra with this dimension vector is isomorphic to~$\S^{a,b}$. It follows, that the function $K_{0,b,a-b,\binom{b+1}{2},b(a-b)}$ is the generating one for the Thom polynomial of this singularity.

\subsection{Nets of conics}

The Thom polynomials for the singularities discussed in this section above have been computed in one or another form by different authors before the method of this paper was developed. Now we present some computations in those cases for which the previous methods were not efficient enough.

The `net of conics' singularity is the one with the local algebra $\frac{\fm_{x,y,z}}{(S_1,S_2,S_3)+\fm^4}$ where $S_1$, $S_2$, and $S_3$ are generic quadratic forms in $x,y,z$.
The following canonical choice for these forms is suggested in~\cite{Wall}:
$$C_\g=\frac{\fm_{x,y,z}}{(S_1,S_2,S_3)+\fm^4}:\qquad S_1=y^2+2xz,
 \quad S_2=2yz,\quad S_3=x^2+6xz+\g z^2.$$
where $\g\in\C$ is a parameter (module). The algebras corresponding to different values of~$\g$ are non-isomorphic. This is the simplest example of a continuous family of non-isomorphic singularities. One of the ways to explain the appearance of a module is the following.
The forms $S_i$ span a three-dimensional subspace in $\Sym^2\C^3$. The intersection of this subspace with the discriminant consisting of degenerate conics defines a cubic projective plane curve, that is, an elliptic curve, whose analytic type is an obvious invariant of the net of conics. For the singularity~$C_\g$ the modular invariant of the discriminantal curve is
$$j=-27\frac{(\g-1)^2}{\g(\g-9)^2}.$$
If $\g\ne0,9$, then the discriminantal curve is smooth, while for $\g=0$ or for $\g=9$ it attains a singularity (a double point).

The algebra~$C_\g$ is $6$-dimensional. Its square defines a filtration with the dimension vector~$(3,3)$. The algebra structure is determined by a linear map
$$\psi:\Sym^2N^{(1)}\to N_2,\qquad N^{(1)}=N_1/N_2,
 \quad\dim N^{(1)}=\dim N_2=3.$$

We obtain that the Thom polynomial of the singularity $C_\g$ (for a fixed value of~$\g$) is given by the generating series
$$C_\g:\qquad K_{3,3}\cdot Z_\g,$$
where $Z_\g$ is the multidegree of the hypersurface ${\mathcal Z}_\g\subset \Alg(3,3)$ consisting of algebras with fixed value of the parameter~$\g$. It is clear that neither~$Z_\g$ nor the Thom polynomial of~$C_\g$ depend on~$\g$. Therefore, it is sufficient to compute them for any distinguished parameter value, say, for $\g=0$ or for $\g=9$.

For $\g=9$ the singularity of the discriminantal curve is attained at the point corresponding to a quadratic form of rank~$1$. It follows that the algebra $C_9$ is characterized by the property that $N^{(1)}$ has a one-dimensional subspace $L$ such that $L\cdot L=0$. The multidegree of the subvariety ${\mathcal Z}_9$ is computed as
$$Z_9=\int_{PN^{(1)}}c_{3}(\Hom(\cO(-2),N_2))$$

For $\g=0$ the singularity of the discriminantal curve is attained at a form of rank~$2$. The simplest way to characterize the locus ${\mathcal Z}_0$ is to apply duality. Let $K\subset \Sym^2N^{(1)}$ be the kernel of the structure map~$\psi$ of the algebra~$C_\g$. By genericity assumption, $\dim K=3$ and $K^\vee\simeq \Sym^2N^{(1)\vee}/N_2^\vee$. Then we have also the restriction map
$$\psi^*:\Sym^2N^{(1)\vee}\to K^\vee.$$
This map can be treated as the structure homomorphism of the algebra $C_{\g^*}$ for another value of the parameter. This defines an involution on the parameter space and one can check that $\g^*=9-\g$. In particular, if $\g=0$ then $\g^*=9$ and by above arguments we get
$$Z_0=\int_{PN^{(1)\vee}}c_{3}(\Hom(\cO(-2),\Sym^2N^{(1)\vee}/N_2^\vee))$$
In both approaches we get
$$Z_0=Z_9=4(c_1(N_2)-2c_1(N^{(1)}))
 =4(t_4+t_5+t_6-2t_1-2t_2-2t_3).$$
We obtained, therefore, a generating series for the net of conics singularity in the form
$$C_\g:\qquad K_{3,3}\cdot 4(t_4+t_5+t_6-2t_1-2t_2-2t_3).$$

The considerations above show that the algebras $C_9$ and $C_0$ posses a reach additional structures. This structures can be used to define more delicate filtrations that lead to more efficient generating series. One can see, for example, that the algebra $C_9$ is generic for the dimension vector~$(0,2,1,0,3)$ or $(0,0,0,1,1,1,0,0,1,1,1)$, and the algebra $C_0$ is generic for the dimension vector $(0,0,1,2,0,1,0,2)$ or $(0, 0, 0, 0, 1, 1, 1, 0, 0, 1, 0, 1, 0, 1)$. This leads to several more possible generating functions for the Thom polynomial of the singularity~$C_\g$:
$$C_\g:\qquad
\begin{array}{l@{\qquad}l}
 K_{0,2,1,0,3},&\frac12K_{0,0,0,1,1,1,0,0,1,1,1},\\
 K_{0,0,1,2,0,1,0,2},&\frac12 K_{0,0,0,0,1,1,1,0,0,1,0,1,0,1}.
\end{array}$$
The factors $\frac12$ in two cases are due to the fact that the corresponding filtrations are not natural in a sense of Definition~\ref{fnat}: they depend on a choice of one of the two branches of the discriminantal curve at its singular point.

Applying Corollary~\ref{cor1} to any of these functions we obtain an expansion for the Thom polynomial
\begin{align*}
\Tp_{C_\g}&=
  4\Delta_{\ell+3,\ell+3,\ell+3,\ell+1,\ell,\ell}
 +8\Delta_{\ell+4,\ell+3,\ell+3,\ell+3,\ell,\ell}
\\&\qquad
 +18\Delta_{\ell+4,\ell+3,\ell+3,\ell+1,\ell,\ell-1}
 +32\Delta_{\ell+4,\ell+4,\ell+3,\ell,\ell,\ell-1}
 +40\Delta_{\ell+4,\ell+4,\ell+3,\ell+1,\ell-1,\ell-1}
\\&\qquad
 +80\Delta_{\ell+4,\ell+4,\ell+4,\ell,\ell-1,\ell-1}
 +32\Delta_{\ell+5,\ell+3,\ell+3,\ell,\ell,\ell-1}
 +20\Delta_{\ell+5,\ell+3,\ell+3,\ell+1,\ell-1,\ell-1}
\\&\qquad
 +120\Delta_{\ell+5,\ell+4,\ell+3,\ell,\ell-1,\ell-1}
 +160\Delta_{\ell+5,\ell+3,\ell+3,\ell-1,\ell-1,\ell-1}
 +80\Delta_{\ell+5,\ell+5,\ell+3,\ell-1,\ell-1,\ell-1}
\\&\qquad
 +40\Delta_{\ell+6,\ell+3,\ell+3,\ell,\ell-1,\ell-1}
 +112\Delta_{\ell+6,\ell+4,\ell+3,\ell-1,\ell-1,\ell-1}
 +16\Delta_{\ell+7,\ell+3,\ell+3,\ell-1,\ell-1,\ell-1}+\dots
\end{align*}
where the dots denote the terms vanishing for $\ell>1$. For the case $\ell=0$ this Thom polynomial is computed in~\cite{FRloc}. It corresponds to the first two terms of the above expansion.

\medskip
We assumed in the above consideration that the parameter~$\g$ of the net of conics is fixed. But one can also consider the singularity type~$C_*$ defined as the union of the singularities $C_\g$ for all~$\g$. Equivalently, the singularity~$C_*$ is identified as a generic one whose nilpotent algebra admits the dimension vector~$(3,3)$. The Thom polynomial of the singularity~$C_*$ is determined by the following generating series
$$C_*:\qquad K_{3,3}.$$

\subsection{Singularities determined by the third-order decomposition}

In this section we explain the computation of Thom polynomials for the Thom-Boardman singularities
$$\S^{2,1,1}=\frac{\fm_{x,y}}{(y^2)+\fm^4},\qquad
 \S^{2,2,1}=\frac{\fm_{x,y}}{(y^3)+\fm^4},\qquad
 \S^{2,2,2}=\frac{\fm_{x,y}}{\fm^4}.$$
These are probably the most interesting examples manifesting the efficiency of the method introduced in this paper. With the previous approaches the problem of finding Thom polynomials for these singularities was considered as a very difficult one or even hopeless. With our approach the answer is obtained almost immediately, with no complicated computations. Consider the gradings on these algebras corresponding to the following choice of the weights of the variables:
\begin{align*}
\S^{2,1,1}:&\quad \deg x=3,\quad \deg y=5,&&
   \ud=(0,0,1,0,1,1,0,1,1,0,1),\\
\S^{2,2,1}:&\quad \deg x=3,\quad \deg y=4,&&
   \ud=(0,0,1,1,0,1,1,1,1,1,1),\\
\S^{2,2,2}:&\quad \deg x=\deg y=1,&&
   \ud=(2,3,4),\\
\S^{2,2,2}:&\quad \deg x=4,\quad \deg y=5,&&
   \ud=(0,0,0,1,1,0,0,1,1,1,0,1,1,1,1).
\end{align*}

\th Proposition. For each of the listed cases, the subvariety in $\Alg(\ud)$ of associative algebras is a complete intersection with a unique component, and a generic algebra from this variety is isomorphic to the corresponding algebra $\S^{2,i,j}$.

The fact that the collection of the associativity equations~\eqref{eqass} form a regular sequence for each of the discussed cases follows from their explicit form. This proofs the first assertion of the proposition. The second assertion for the singularities $\S^{2,2,1}$ and $\S^{2,2,2}$ are particular cases of Proposition~\ref{propquasi}. That proposition does not cover the case of the singularity~$\S^{2,1,1}$ and the dimension vector~$(0,0,1,0,1,1,0,1,1,0,1)$ since the list of monomials of degree~$\le11$ in~$x$ and~$y$ with~$\deg x=3$, $\deg y=5$ contains also $y^2$ which is zero in the algebra~$\S^{2,1,1}$. The arguments in the proof of Proposition~\ref{propquasi} show that a generic associative algebra with this dimension vector can be written in the form
 $$N=\frac{\<x,y,x^2,xy,x^3,y^2,x^2y\>}{(y^2+a\,x^2y)}$$
for some $a\in\C$. Applying the change of coordinates $y=\tilde y-\frac12a\,x^2$ we can bring the relation $y^2+a\,x^2y=0$ to the form $\tilde y^2=0$. It proves that the algebra~$N$ is isomorphic to~$\S^{2,1,1}$ for any~$a$.\QED

The proposition implies the formulas for the generating functions of Thom polynomials for these singularities from the table of Sect.~\ref{sectab}. An extra factor $t_1$ in the case of the dimension vector $(0,0,0,1,1,0,0,1,1,1,0,1,1,1,1)$ for the singularity~$\S^{2,2,2}$ is due to the fact that a filtration with this dimension is not natural, its choice is determined by a choice of the direction of the $y$-axis in the plane of coordinates~$x$ and~$y$. Possible filtrations form the projective line~$\C P^1$ and the normalizing monomial~$t_1$ is chosen by the requirement $\int_{\C P^1}t_1=1$ in accordance with the Definition~\ref{fnat}.

One can check the obtained Thom polynomials by restricting them to some special cases. For example, setting $\ell=-1$, $1+c_1+c_2+\dots=\frac{1+\b}{(1+\a_1)(1+\a_2)}$ in the Thom polynomial for $\S^{2,2,2}$ we get its expected value for the case of the $2$-dimensional source and $1$-dimensional target of the mapping,
$$\Tp_{\S^{2,2,2}}\bigm|_{m=2,n=1}=c_{\rm top}(J^3_0(\C^2,\C))=
  \prod_{1\le i+j\le 3}(\b-i\a_1-j\a_2).$$

\section{Final remarks and open questions}\label{sec9}

A complete classification of singularities is so wild and irregular that there is no hope to obtain the Thom polynomial for any given in advance singularity in a closed form. I have even doubts whether the Thom polynomial of the Morin singularity, $A_\mu$, will ever be computed for all~$\mu$, without speaking about generic Thom-Boardman singularities~$\S^{j_1,\dots,j_r}$. On can extend the list of singularities with known Thom polynomials, say, by classifying singularities with nilpotent algebras of small dimension, or those determined by small order of jets. One should notice, however, that the new classes of singularities have a relatively big codimension. Any explicit computations for them would require considerable computer resources which makes them less attractive from the viewpoint of applications.

One of the possible extensions of the presented approach is to the multisingularity theory~\cite{Kaz03}. In this theory, the class of a given multisingularity type of a mapping is expressed in terms of the so called residual polynomials in the relative Chern classes of the mapping. First experiments show that the residual polynomials of multisingularities experience a stabilization with the growth of the relative dimension~$\ell$ similar to that for Thom polynomials of monosingularities, but there is to general theorem explaining this stabilization.

The non-uniqueness of the generating series brings up a number of interesting questions on the algebra of generating functions. Let us call two rational functions in~$t_1,\dots,t_\mu$ \emph{equivalent} if they lead to equal expansions in the Chern classes. Is there a simple way to determine whether two rational functions are equivalent? What is the smallest number of linear factors in the denominator of rational functions equivalent to the given one? These questions will be addressed in further studies.


\end{document}